\documentclass[reqno,11pt]{article}
\usepackage{amsmath,amsfonts,amssymb,euscript,theorem}

\numberwithin{equation}{section}
\newtheorem{theorem}{Theorem}[section]

\newtheorem{define}[theorem]{Definition}
\newtheorem{corollary}[theorem]{Corollary}
{\theorembodyfont{\rmfamily} \newtheorem{remark}[theorem]{Remark} }
\newtheorem{lemma}[theorem]{Lemma}
\newtheorem{proposition}[theorem]{Proposition}

\newcommand{\bke}[1]{\left( #1 \right)}

\newcommand{\bket}[1]{\left\{ #1 \right\}}
\newcommand{\norm}[1]{\left\Vert #1 \right\Vert}
\newcommand{\abs}[1]{\left| #1 \right|}
\newcommand{\bka}[1]{\left\langle #1 \right \rangle}

\newcommand{\al}{\alpha}
\newcommand{\e}{\varepsilon}
\newcommand{\de}{\delta}
\newcommand{\ga}{\gamma}

\newcommand{\si}{\sigma}
\renewcommand{\th}{\theta}
\newcommand{\R}{\mathbb{R}}
\newcommand{\C}{\mathbb{C}}
\newcommand{\Z}{\mathbb{Z}}
\newcommand{\bS}{\mathbb{S}}
\renewcommand{\Re}{\mathop{\mathrm{Re}}}
\renewcommand{\Im}{\mathop{\mathrm{Im}}}
\newcommand{\pd}{\partial}
\newcommand{\donothing}[1]{{}}
\newcommand{\mat}[1]{\begin{bmatrix}#1 \end{bmatrix}}

\newcommand{\lec}{{\ \lesssim \ }}

\newcommand{\sE}{\EuScript{E}}
\newcommand{\sO}{\EuScript{O}}

\newcommand{\Emin}{\sE_\text{min}}

\newcommand{\bv}{{\bf v}}
\newcommand{\bu}{{\bf u}}
\newcommand{\bh}{{\bf h}}
\newcommand{\bxi}{\boldsymbol{\xi}}
\newcommand{\bta}{\boldsymbol{\eta}}
\newcommand{\bV}{{\bf V}}
\newcommand{\bphi}{\boldsymbol{\phi}}

\newcommand{\p}{\partial}

\newcommand{\hi}{\widehat {\pmb {\imath}}}
\newcommand{\hj}{\widehat {\pmb {\jmath}}}
\newcommand{\hk}{\widehat {\mathbf{k}}}
\newcommand{\be}{\mathbf{e}}
\newcommand{\BE}{ \mathbf{\hat E}}
\newcommand{\BV}{ \mathbf{V}}

\newcommand{\Tmax}{T_{\rm{max}}}
\newcommand{\nrm}[1]{\| #1 \|}

\newenvironment{proof}{{\bf Proof.}}{\hfill\fbox{}\par\vspace{.2cm}}

\begin{document}

\title{Asymptotic stability of harmonic maps \\ under the Schr\"odinger
flow 
}
\author{{}\\[0mm]
Stephen Gustafson \quad\quad Kyungkeun Kang \quad\quad Tai-Peng Tsai
}

\date{}
\maketitle


\renewcommand{\S}{\mathbb{S}}

\begin{abstract}
  For Schr\"odinger maps from
  $\R^2\times\R^+$ to the $2$-sphere $\S^2$, it is not known if
  finite energy solutions can form singularities (``blowup'') in finite time.
  We consider equivariant solutions with energy near the energy
  of the two-parameter family of equivariant harmonic maps.
  We prove that if the topological degree of the map is at least
  four, blowup does {\it not} occur, and global solutions converge
  (in a dispersive sense -- i.e. scatter) to a fixed harmonic map
  as time tends to infinity.
  The proof uses, among other things, a time-dependent
  splitting of the solution, the ``generalized Hasimoto transform",
  and Strichartz (dispersive) estimates for a certain
  two space-dimensional linear Schr\"odinger equation whose potential
  has critical power spatial singularity and decay.
  Along the way, we establish an energy-space local well-posedness
  result for which the existence time is determined by the length-scale
  of a nearby harmonic map.
\end{abstract}


\tableofcontents

\section{Introduction and main results}
\label{sec:intro}

The {\it Schr\"odinger flow} for maps from $\R^n$ to $\S^2$
(also known as the {\it Schr\"odinger map}, and, in ferromagnetism, as the
{\it Heisenberg model} or {\it Landau-Lifshitz equation})
is given by the equation
\begin{equation}
  \frac{\p {\bf u}}{\p t} = {\bf u} \times \Delta {\bf u},
  \quad\quad  {\bf u}(x,0)={\bf u}_0(x).
\label{eq:sm}
\end{equation}
Here ${\bf u} = {\bf u}(x,t)$ is the unknown map from
$\R^n \times \R^+$ to the $2$-sphere
\[
  \S^2 := \{ {\bf u} \in \R^3 \;\; | \;\; |{\bf u}|=1 \} \subset \R^3,
\]
$\Delta$ denotes the Laplacian in $\R^n$, and
$\times$ denotes the cross product of vectors in $\R^3$.
A somewhat more geometric way of writing Equation~\eqref{eq:sm} is
\begin{equation}
  \frac{\p {\bf u}}{\p t} = J P \Delta {\bf u}
\end{equation}
where $P = P^{{\bf u}}$ denotes the orthogonal projection from
$\R^3$ onto the tangent plane
\[
  T_{\bf u} \S^2 := \{ \bxi \in \R^3 \; | \; \bxi \cdot {\bf u} = 0 \}
\]
to $\S^2$ at ${\bf u}$ (so that
$P \Delta {\bf u} = \Delta {\bf u} + |\nabla {\bf u}|^2 {\bf u})$, and
\[
  J = J^{\bf u} := {\bf u} \times
\]
is a rotation through $\pi/2$ on the tangent plane $T_{\bf u} \S^2$.

On one hand, Equation~\eqref{eq:sm} is a borderline case of the
{\it Landau-Lifshitz-Gilbert} equations which model dynamics
in isotropic ferromagnets (including dissipation):
\begin{equation}
  \frac{\p {\bf u}}{\p t} =
  a P \Delta {\bf u} + b J P \Delta {\bf u},
  \quad\quad a \ge 0
\end{equation}
(see, eg., \cite{KIK}).
The Schr\"odinger flow corresponds to the case $a=0$. The case $b=0$
is the well-studied harmonic map heat flow, for which some
finite-energy solutions do blow up in finite time (\cite{cdy}).

On the other hand, Equation~\eqref{eq:sm} is a particular case of
the Schr\"odinger flow for maps from a Riemannian manifold into a
K\"ahler manifold (see, eg., \cite{DW98,TU,GS,D02}). We will
consider only the case of maps $: \R^2 \times \R^+ \to \S^2$ in
this paper.

We refer the reader to our previous paper~\cite{GKT} for
more detailed background on~\eqref{eq:sm} (and further references),
limiting the discussion
here to a list of a few basic facts we need in order to state our results.

\begin{itemize}
\item {\bf Energy conservation.}
Equation \eqref{eq:sm} formally conserves the {\it energy}
\begin{equation}
  \sE({\bf u}) := \frac 12 \int_{\R^n} |\nabla {\bf u}|^2 \, dx
  = \frac 12 \int_{\R^n} \sum_{j=1}^n \sum _{k=1}^3
  \left| \frac{\p u_k}{\p x_j} \right|^2 \, dx.
\end{equation}
The space dimension $n=2$ is critical in the
sense that $\sE({\bf u})$ is invariant under scaling.  In general,
\begin{equation}
  \sE({\bf u}(\cdot)) = s^{2-n} \sE({\bf u}(\cdot/s))
\end{equation}
for $s > 0$.

\item {\bf Equivariant maps.}
Fix $m \in \Z$ a non-zero integer.
By an {\it $m$-equivariant map} ${\bf u}:\R^2 \to \S^2 \subset \R^3$,
we mean a map of the form
\begin{equation}
 {\bf u}(r,\th) =  e^{m \th R} \, {\bf v}(r)
\end{equation}
where $(r,\th)$ are polar coordinates on $\R^2$,
${\bf v} : [0, \infty) \to \S^2$, and $R$ is the
matrix generating rotations around the $u_3$-axis:
\begin{equation}
R=\begin{bmatrix}0&-1&0\\1&0&0\\0&0&0\end{bmatrix},
\quad
e^{\al R}=\begin{bmatrix}\cos\al&-\sin\al&0\\
\sin\al&\cos\al&0\\0&0&1\end{bmatrix}.
\end{equation}
Radial maps arise as the case $m=0$.
The class of $m$-equivariant maps is formally
preserved by the Schr\"odinger flow.

\item {\bf Topological lower bound on energy.}
If ${\bf u}$ is $m$-equivariant, we have
$|\nabla {\bf u}|^2 = |\p {\bf u}/\p r|^2
+ r^{-2} |\p {\bf u}/\p \th|^2
= |\p {\bf v}/\p r|^2 + (m^2/r^2)|R{\bf v}|^2$ and so
\begin{equation}
  \sE({\bf u}) = \pi\int_0^\infty
  \bke{ \left| \frac{\p {\bf v}}{\p r} \right|^2
  + \frac {m^2}{r^2} (v_1^2+v_2^2)} \, r dr .
\end{equation}
If $\sE({\bf u})<\infty$, then ${\bf v}(r)$ is continuous,
and the limits $\lim_{r \to 0} {\bf v}(r)$ and
$\lim _{r \to \infty} {\bf v}(r)$ exist
(see \cite{GKT}), and so we must have
${\bf v}(0),{\bf v}(\infty) = \pm \hk$, where $\hk = (0,0,1)^T$.
Without loss of generality we fix ${\bf v}(0)=-\hk$.
The two cases ${\bf v}(\infty) = \pm \hk$ then
correspond to different topological classes of maps.
We denote by $\Sigma_m$ the class of $m$-equivariant maps with
${\bf v}(\infty) = \hk$:
\begin{equation}
  \Sigma_m = \bket{{\bf u}: \R^2 \to \S^2 \;\; | \;\;
  {\bf u} = e^{m \th R}{\bf v}(r), \;\;
  \sE({\bf u})<\infty, \;\; {\bf v}(0) = -\hk, \;\;
  {\bf v}(\infty) = \hk}.
\end{equation}
For ${\bf u} \in \Sigma_m$, the energy $\sE({\bf u})$ can be rewritten:
\begin{equation}
\label{eq:energy2}
  \sE({\bf u}) = \pi \int_0^\infty
  \bke{ \left| \frac{\p {\bf v}}{\p r} \right|^2 +
  \frac {m^2}{r^2} |J^{\bf v}  R {\bf v}|^2} \, r dr
  = \pi \int_0^\infty \left| \frac{\p {\bf v}}{\p r}
  - \frac{|m|}{r} J^{\bf v} R {\bf v} \right|^2 \, r dr + \Emin
\end{equation}
(recall $J^{\bf v} := {\bf v} \times$) with
\begin{equation}
  \Emin = 2\pi \int_0^\infty {\bf v}_r \cdot
  \frac{|m|}{r} J^{\bf v} R {\bf v} \, r dr
  = 2 \pi |m| \int_0^\infty (v_3)_r dr = 4 \pi |m|.
\end{equation}
Thus for ${\bf u} \in \Sigma_m$, there is a
non-trivial lower bound for the energy:
\begin{equation}
\label{eq:tlb}
  {\bf u} \in \Sigma_m \quad \implies \quad
  \sE({\bf u}) \geq 4 \pi |m|.
\end{equation}
(In general one has $\sE({\bf u}) \geq 4\pi|deg|$
where $deg$ is the {\it degree} of the map ${\bf u}$,
considered as a map from $\S^2$
to itself (defined, for example, by integrating the pullback
by ${\bf u}$ of the volume form on $\S^2$).)

\item {\bf Harmonic maps.}
For a map ${\bf u} \in \Sigma_m$,
the topological lower bound~\eqref{eq:tlb} is saturated if and only if
\begin{equation}
\label{eq:hmeq}
  \frac{\p {\bf v}}{\p r} = \frac{|m|}{r} J^{\bf v} R {\bf v},
\end{equation}
and the minimal energy is attained
(i.e. \eqref{eq:hmeq} is satisfied) precisely at
the two-parameter family of harmonic maps
\begin{equation}
\label{eq:hm}
  \sO_m := \bket{e^{(m \th + \al)R} \bh(r/s) \;\;
  | \;\; s > 0, \; \al \in \R}
\end{equation}
where
\begin{equation}
\label{eq:harmonic}
  \bh(r) = \left( \begin{array}{c} h_1(r) \\ 0 \\ h_3(r)
  \end{array} \right), \quad
  h_1(r)=\frac 2{r^{|m|} + r^{-|m|}}, \quad
  h_3(r)= \frac {r^{|m|} - r^{-|m|}}{r^{|m|} + r^{-|m|}}.
\end{equation}
The rotation parameter $\alpha$ is determined
only up to shifts of $2 \pi$ (i.e. really $\alpha \in \S^1$).
The fact that $\bh(r)$ satisfies~\eqref{eq:hmeq} means
\begin{equation}
  (h_1)_r = -\frac m r h_1 h_3, \quad (h_3)_r = \frac m r h_1^2.
\end{equation}
Note that $\sO_m$ is just the orbit of the harmonic map
$e^{m \theta R} \bh(r)$ under the symmetries of the energy
$\sE$ which preserve equivariance: scaling and rotation.
Explicitly, the maps in $\sO_m$ are of the form
\begin{equation}
  u(r, \theta) = \left(
  \begin{array}{c} \cos (m \th + \al)h_1( r/s) \\
  \sin(m \th + \al) h_1(r/s) \\
   h_3( r/s) \end{array} \right).
\end{equation}
Of course, these harmonic maps are each static
solutions of the Schr\"odinger flow \eqref{eq:sm}.
In fact, it is not hard to show they are the only 
$m$-equivariant static solutions (though this fact plays
no role in our analysis).

\item {\bf The ``orbital stability'' of $\sO_m$.}
We recall the main result of~\cite{GKT}:
\begin{theorem}\cite{GKT}
\label{thm:gkt}
There exist $\delta > 0$ and $C_1, C_2 > 0$ such that if
${\bf u} \in C([0,T) ; \dot{H}^2 \cap \Sigma_m)$ is
a solution of the Schr\"odinger flow \eqref{eq:sm}
conserving energy, and satisfying
\[
  \delta_1^2 := \sE({\bf u}_0) - 4\pi |m| < \delta^2,
\]
then there exist $s(t)\in{\mathcal C}([0,T);(0,\infty))$ and
$\alpha(t)\in{\mathcal C}([0,T);\R)$ so that
\begin{equation}
\label{eq:closeness}
  \norm{{\bf u}(x,t)-e^{(m\theta+\alpha(t))R}\bh(r/s(t))}_{\dot{H}_1(\R^2)}
  \leq C_1 \delta_1, \quad \forall t \in [0,T).
\end{equation}
Moreover, $s(t) > C_2/\|{\bf u}(t)\|_{\dot{H}_2(\R^2)}$.
Furthermore, if $T < \infty$ is the maximal time of existence
for ${\bf u}$ in $\dot{H}^2$
(i.e. if $\lim_{t \to T^-} \| {\bf u}(t) \|_{\dot H^2(\R^2)} = \infty$), then
\begin{equation}
\label{eq:concentrate}
  \liminf_{t\to T^{-}} s(t) = 0.
\end{equation}
\end{theorem}
This theorem can be viewed, on one hand, as an
{\it orbital stability} result for the family $\sO_m$
of harmonic maps (at least up to the possible blowup
time), and on the other hand as a characterization
of blowup for energy near $\Emin$:
solutions blowup if and only if the ``length-scale''
$s(t)$ goes to zero. Here $s(t)$ (and the rotation
angle $\alpha(t)$) are determined simply by finding,
at each time $t$, the harmonic map which is
$\dot{H}^1$-closest to $\bu(t)$. More precisely,
a continuous map
\begin{equation}
\label{eq:sdef}
\begin{split}
  \{ \; {\bf u} \in \Sigma_m \; | \; \sE({\bf u}) < 4 \pi |m| + \delta^2 \; \}
  & \to \R^+ \times (\R \mod 2\pi) \\
  {\bf u} & \mapsto ( \; s(\bf{u}), \; \alpha({\bf u}) \; )
\end{split}
\end{equation}
is constructed in \cite{GKT},
which, for $m$-equivariant maps with energy close to
$4 \pi |m|$, identifies the unique $\dot H^1$-closest
harmonic map:
\begin{equation}
\label{eq:sdef2}
  \| {\bf u} - e^{[m \theta + \alpha({\bf u})]R} h(r/s({\bf u})) \|_{\dot H^1}
  = \min_{s \in \R^+, \alpha \in \R}
  \| {\bf u} - e^{[m \theta + \alpha]R} h(r/s) \|_{\dot H^1}.
\end{equation}
Then we set $s(t) := s(\bu(t))$. 
\end{itemize}

In this paper, we continue our study of the
Schr\"odinger flow for equivariant maps with energy close to the harmonic
map energy. 
We begin with an energy-space local well-posedness theorem for such maps.
It is worth remarking that despite a great deal of recent work on the local
well-posedness problem in two space dimensions
(\cite{SSB,DW01,HM,B,IK}; see also \cite{NSU,KK,Kenig}
for the ``modified Schr\"odinger map'' case),
there is no general result for energy space initial data. 
For our special class of data, however, we do have such a result.
Before stating it, let us first make precise the sense
in which our energy-space solution solves the Schr\"odinger map
problem:
\begin{define}[Weak solutions]
\label{def:sol}
Let $Z := \{ \; \bu : \R^n \to \bS^2, \; D \bu \in L^2 \; \}$
be the energy space.
We say $\bu(x,t)$ is a {\it weak solution}
of the Schr\"odinger flow~\eqref{eq:sm} 
on the time interval $I = [0,T]$, with initial data
$\bu_0 \in Z$, if
\begin{enumerate}
\item 
$\bu \in L^\infty(I;Z) \cap C_{weak}([0,T];Z)$ 
\item $\bu(0) = \bu_0$
\item 
$\iint_{\R^n \times I} \{ \bu \cdot\bphi_t - \bu \times \pd_j \bu \cdot \pd_j \bphi \}
dx dt = 0$
for all $\bphi \in C^1_c(I \times \R^n; \R^3)$.
\end{enumerate}
\end{define}
\begin{remark}
It is not strictly necessary
to require that $D \bu$ be weakly continuous in $t$
(in property $1$ above). 
The weak form of the equation (property $3$)
implies $\bu_t \in L^\infty([0,T];H^{-1})$,
and so, after redefinition on a set of
time measure zero, $\bu \in Lip([0,T];H^{-1})$
and $D \bu \in Lip([0,T];H^{-2})$.
Since we also have $D \bu \in L^\infty([0,T];L^2)$,
we can prove $D \bu \in C_{weak}([0,T]; L^2)$.
\end{remark}

We have
\begin{theorem}[Local wellposedness]
\label{lwp:T1}
Let $|m| \ge 1$.  There exist $\delta > 0$ and $\si>0$ such that the
following hold.  Suppose $\bu_0 \in \Sigma_m$ and 
$\sE(\bu_0) = 4 \pi m + \de_0^2$, $\de_0 \in (0, \de]$.  
Let $s_0 := s(\bu_0)$, as defined
in~\eqref{eq:sdef}-\eqref{eq:sdef2}. 
Then there is a unique weak solution $\bu(t)$ of \eqref{eq:sm}
\[
  \bu(t) \in C(I; \Sigma_m), \quad I=[0,\si s_0^2].
\]
Moreover, $\sE(\bu(t)) = \sE(\bu_0)$ for $ t \in I$. If, furthermore, $\bu_0
\in \dot H^2$, then $\bu(t) \in C(I; \Sigma_m \cap \dot H^2)$.
Suppose $\bu^n_0\to \bu_0$ in $\Sigma_m$
and let $\bu^n$ denote the corresponding solutions of \eqref{eq:sm},
then $\bu^n\to \bu$ in $C(I,\Sigma_m)$.
\end{theorem}

It is worth emphasizing that the existence time furnished
by this theorem depends not on the energy
$\|\bu_0\|^2_{\dot H^1}$ of the initial data
(reflecting the energy-space critical nature of the equation in
dimension $n=2$), but rather on $s(\bu_0)$, the length scale of the
$\dot H^1$-nearest harmonic map.

There are at least two ways to define blow-up for these
solutions. Suppose $\bu(t) \in C([0,T), \Sigma_m \cap \dot H^k)$,
$0<T<\infty$ with $k=1$ or $2$.  If
$k=1$, we say $\bu(t)$ blows up at $t=T$ if $\lim_{t \to T^-} \bu(t)$
does not exist in $\dot H^1$.
If $k=2$, we say $\bu(t)$ blows up at $t=T$
if $\limsup_{t \to T-} \norm{\bu(t)}_{\dot H^2} = \infty$.

For $\bu_0 \in \Sigma_m\cap \dot H^k$, $k=1,2$, denote by $\Tmax^k$ the
maximal time such that there is a unique solution $\bu(t) \in
C([0,\Tmax^k); \Sigma_m\cap \dot H^k)$.

\begin{corollary}
\label{lwp:T1b}
Under the same assumptions as in Theorem \ref{lwp:T1}, suppose the
solution $\bu(t) \in C([0,T),\Sigma_m \cap \dot{H}^k)$, $k=1$ or $2$,
and $T < \infty$.

\begin{itemize}
\item[(i)] (Blowup alternative) $\bu(t)$ blows up at time $T$ 
(i.e. $T = T^1_{max}$) iff $\liminf_{t \to T-} s(\bu(t)) = 0$.  
In this case, $s(\bu(t)) \le C \sqrt{T-t}$, and
if $k=2$, $T = T^1_{max} = T^2_{max}$ with
$\norm{\bu(t)}_{\dot H^2} \ge C(T-t)^{-1/2}$.

\item[(ii)] (Lower bound for $T_{\rm{max}} := T^1_{max}$) We have 
$\Tmax \ge \si [s(\bu_0)]^2$ 
(here $\sigma$ is the constant from Theorem~\ref{lwp:T1}).
\end{itemize}
\end{corollary}
Corollary~\ref{lwp:T1b} (i) improves Theorem~\ref{thm:gkt} by
giving explicit bounds. 

We also have $\dot H^1$ local wellposedness for the small energy
equivariant case considered in \cite{CSU}. Since the energy is conserved, local
wellposedness implies global wellposedness.

\begin{theorem}[Small energy local wellposedness]
\label{lwp:T4}
Let $|m| \ge 1$.  There exist $\delta > 0$ and $\si>0$ such that the
following hold.  Suppose $\bu_0 = e^{m \th R} v_0(r)$ and $\sE(\bu_0) \le
\de^2$, then there is a unique weak solution $\bu(t,r,\th)= e^{m \th R}
\bv(t,r)$ of \eqref{eq:sm} so that $\bu(t) \in C([0,\si]; \dot H^1)$.
Moreover, $\sE(\bu(t)) = \sE(\bu_0)$ for $t \in [0,\si]$.
Suppose $\bu^n_0$ are equivariant, $\bu^n_0\to \bu_0$ in $\dot H^1$ and
let $\bu^n$ denote the corresponding solutions of \eqref{eq:sm}, then
$\bu^n\to \bu$ in $C([0,\si],\dot H^1)$.
\end{theorem}
Note that this result does {\it not} cover the radial 
case ($m=0$).

The question of whether singularities can form in the
Schr\"odinger flow is open. So far, it has only been shown that
they cannot form for {\it small energy} radial or equivariant
solutions (\cite{CSU}).
Our Theorem~\ref{thm:gkt} above leaves open the
question of whether finite-time blowup can occur
for maps in $\Sigma_m$ with energies near $\Emin = 4\pi|m|$. 
The main result of this paper
shows that when $|m| \geq 4$, it does {\it not}.
Moreover, we show that these solutions converge
(in a dispersive sense) to specific harmonic maps
as $t \to \infty$.  Here is the main result:
\begin{theorem}[Main result]
\label{thm:gkt2}
Let $|m| \geq 4$.
Let $(r,p)$ satisfy $2 < r \leq \infty$, $2 \leq p < \infty$,
with $1/r + 1/p = 1/2$.
There exist positive constants 
$\delta$, $C$, and $C_p$, such that if
${\bf u}_0 \in \Sigma_m$ satisfies
\[
  \delta_1^2 := \sE({\bf u}_0) - 4\pi |m| < \delta^2,
\]
then for the corresponding solution
${\bf u}(t)$ of the Schr\"odinger flow
(guaranteed by Theorem~\ref{lwp:T1}),
\begin{enumerate}
\item
there is no finite-time blowup: $\Tmax = \infty$
\item
there exist
$s(t)\in{\mathcal C}([0,\infty);(0,\infty))$ and
$\alpha(t)\in{\mathcal C}([0,\infty);\R)$ such that
\begin{equation}
\label{eq:thm1}
  \norm{ \nabla [ {\bf u}(x,t)-e^{(m\theta+\alpha(t))R}\bh(r/s(t)) ] }
  _{(L^\infty_t L^2_x \cap L^r_t L^p_x)(\R^2 \times [0, \infty))}
  \leq C_p \delta_1
\end{equation}
\item
furthermore,
\[
  \left| \frac{s(t)}{s(\bu_0)} - 1 \right| +
  |\alpha(t) - \alpha(\bu_0)| \leq C \delta_1^2
  \quad\quad \mbox{ for all } \; t > 0
\]
and there exist $s_+ > 0$ and $\alpha_+$ with
\begin{equation}
\label{eq:thm2}
  s(t) \to s_+,  \quad
  \alpha(t) \to \alpha_+, \quad
  \mbox{ as } \;\; t \to \infty.
\end{equation}
\end{enumerate}
\end{theorem}

\begin{remark}
\begin{enumerate}
\item 
The $L^\infty_t L^2_x$ (energy space) estimate
in~(\ref{eq:thm1}) already follows from Theorem~\ref{thm:gkt}. The
other space-time estimates in~(\ref{eq:thm1}) further imply
asymptotic {\it convergence} to the family of harmonic maps (at
least, in a time-averaged sense -- the best we can expect without
further assumptions on the initial data). The convergence
results~\eqref{eq:thm1} and~\eqref{eq:thm2} are precisely
what we mean when we say the harmonic maps are 
{\it asymptotically stable} under the
Schr\"odinger flow for $|m| \geq4$.
\item 
Note that for $|m|=1, 2, 3$, the fate of solutions with
energy near $\Emin$ is still an open question. Our restriction
$|m| > 3$ is connected with the slow spatial decay of the harmonic
map component $h_1(r) \sim (const) r^{-|m|}$ as $r \to \infty$.
For a somewhat technical reason, we need $r^2 h_1(r) \in L^2(rdr)$
(see Lemma~\ref{lem:zbyq}), which requires $|m| > 3$. For 
seemingly more fundamental reasons, we need $rh_1(r) \in L^2(rdr)$
(see~\eqref{eq:restriction}), which holds if $|m| > 2$.
\item 
The recent work \cite{RodSter} on the analogous
{\it wave map} problem, imposes the same $|m| \geq 4$
restriction, but proves that {\it blow-up} is possible
in this class, suggesting that singularity formation
is a more delicate question for Schr\"odinger maps
than for wave maps.  
\end{enumerate}
\end{remark}

We end the introduction with a few words about our approach.
One key observation, already used in~\cite{GKT}, is that the
tangent vector field
\[
  {\bf W} := \frac{\p {\bf v}}{\p r} - \frac{|m|}{r}J^{\bf v} R {\bf v}
\]
``measures the deviation of the map $\bu$ from harmonicity''
(this is indicated by~\eqref{eq:hmeq}, for example).
Furthermore, when expressed in an appropriate orthonormal frame,
the coordinates of ${\bf W}$ satisfy a nonlinear Schr\"odinger-type equation
which is suitable for obtaining estimates --
this is the {\it generalized Hasimoto transform} introduced
in \cite{CSU} to study the small energy problem.

In the present work, this nonlinear Schr\"odinger-type PDE
is coupled to a two-dimensional dynamical system
describing the dynamics of the
scaling and rotation parameters $s(t)$ and $\alpha(t)$, a careful choice
of which must be made at each time in order to allow estimation.
This is all done in Section~\ref{sec:system}.

The key to proving convergence of the solution to a harmonic map
is then to obtain dispersive estimates --
in this case Strichartz-type estimates -- for the linear part of
our nonlinear Schr\"odinger equation.
The potential appearing in the corresponding Schr\"odinger operator turns
out to have $const/|x|^2$ behaviour both at the origin, and as
$|x| \to \infty$, which is a ``borderline'' case not treatable
by purely perturbative methods. Fortunately, a recent series
of papers by Burq, Planchon, Stalker, and Tahvildar-Zadeh
(see \cite{BPST1,BPST4}) addresses the problem of obtaining
dispersive estimates when the potential has just this ``critical''
decay rate, provided the potential satisfies a
``repulsivity'' condition (which in particular rule out bound states).
Though their relevant results are for dimension $n \geq 3$,
we are able to adapt their approach to prove the estimates we need
in our two-dimensional setting. This is done in Section~\ref{sec:linear}.

Finally, in Section~\ref{sec:proof}, we prove Theorem~\ref{thm:gkt2}
by applying the linear estimates of Section~\ref{sec:linear}
to the coupled nonlinear system of Section~\ref{sec:system}.

Since the proof of Theorems~\ref{lwp:T1} and~\ref{lwp:T4}, and
Corollary \ref{lwp:T1b} are independent of the rest of the paper,
they are postponed to Section~\ref{sec:lwp}.
Some lemmas are proved in Section~\ref{sec:lemmas}.

\begin{remark}
\begin{enumerate}
\item
From here on, we will assume $m>0$. For $m<0$, simply
make the change of variable $(x_1,x_2,x_3) \to (x_1,-x_2,x_3)$.
\item {\bf Notation:}
throughout the paper, the letter
$C$ is used to denote a generic constant,
the value of which may change from line to line.
Vectors in $\R^3$ appear in boldface, while their
components appear in regular type: for example,
${\bf u} = (u_1, u_2, u_3)$.
\end{enumerate}
\end{remark}


\section{The dynamics near the harmonic maps}
\label{sec:system}

\subsection{Splitting the solution}
\label{sec:splitting}

Let ${\bf u}(x,t) = e^{m\theta R} {\bf v}(r,t) \in \Sigma_m$
be a solution of the Schr\"odinger map equation~\eqref{eq:sm}.
We will write our solution as a harmonic map with time-varying
parameters, plus a perturbation:
\begin{equation}
\label{eq:split}
  {\bf v}(r,t) = e^{\alpha(t) R}
  \left[ \bh(\rho) + \bxi(\rho,t) \right], \quad\quad
  \rho := \frac{r}{s(t)}
\end{equation}
In Section~\ref{sec:orth} we take up the central
question of precisely how to do this splitting (i.e. the choice
of $s(t)$ and $\al(t)$).

It is convenient and natural to single out the component
of the perturbation $\bxi$ which is tangent to $\bS^2$ at $\bh$:
\[
  \bxi(\rho,t) = \bta(\rho,t) + \gamma(\rho,t) \bh(\rho),
  \quad\quad \bta(\rho,t) \in T_{\bh(\rho)} \bS^2,
\]
so that $\bta \cdot \bh \equiv 0$. Thus the original map
$\bu$ is written
\[
  \bu(x,t) = e^{[m\theta + \alpha(t)]R}
  [ (1 + \gamma(\rho,t))\bh(\rho) + \bta(\rho,t) ];
  \quad\quad \rho = \frac{r}{s(t)}, \quad\quad
  \bta(\rho,t) \in T_{\bh(\rho)} \bS^2.
\]
The pointwise constraint $|{\bf v}| \equiv 1$ forces
\[
  1 \equiv |\bh + \bxi|^2 = |(1+\gamma)\bh + \bta|^2
  = (1 + \gamma)^2 + |\bta|^2,
\]
so $\gamma(\rho,t) \leq 0$ and $|\bta(\rho,t)| \leq 1$.
If $|\bxi| \leq 1$, then
\begin{equation}
\label{eq:gamma1}
  \gamma(\rho,t) = +(1 - |\bta(\rho,t)|^2)^{1/2} - 1 \in [-1,0],
\end{equation}

A convenient orthonormal basis of $T_{\bh(\rho)} \bS^2$ is given by
\[
  \hj := \left( \begin{array}{c} 0 \\ 1 \\ 0 \end{array} \right),
  \quad \mbox{ and } \quad J^{\bh(\rho)} \hj =
  \left( \begin{array}{c} -h_3(\rho) \\ 0 \\ h_1(\rho) \end{array} \right),
\]
and we will express tangent vectors like
$\bta \in T_\bh \bS^2$ in this basis via the invertible linear map
\[
\begin{split}
  \bV^\rho : \C &\to T_{\bh(\rho)} \bS^2 \\
  z = z_1 + iz_2 &\mapsto z_1 \hj + z_2 J^{\bh(\rho)} \hj.
\end{split}
\]
So we write
\[
  \bta(\rho,t) = \bV^\rho(z(\rho,t)),
\]
and in this way, the complex function $z(\rho,t)$,
together with a choice of the parameters $s(t)$ and $\alpha(t)$,
gives a full description of the original solution ${\bf u}(x,t)$,
provided $|\bxi| \leq 1$.

From~\eqref{eq:gamma1}, we find
\begin{equation}
\label{eq:gamma2}
  |z| = |\bta| \leq 1/2 \quad \implies \quad
  |\gamma| \lec |z|^2, \quad
  |\gamma_\rho| \lec |z| |z_\rho|.
\end{equation}
These estimates, together with results in~\cite{GKT},
show that if $s$ and $\al$ are chosen appropriately, then
for $\sE(\bu) - 4\pi m$ small,
\[
  \| z \|_X^2 \lec \sE(\bu) - 4\pi m \lec \| z \|_X^2
\]
where
$X := \{ z : [0,\infty) \to \C \; | \;
z_\rho \in L^2(\rho d\rho), \;\; \frac{z}{\rho} \in L^2(\rho d\rho)\}$,
with
\begin{equation}
\label{eq:X}
  \| z \|_X^2 := \int_0^\infty \left\{
  |z_\rho(\rho)|^2 + \frac{|z(\rho)|^2}{\rho^2} \right\} \rho d\rho.
\end{equation}
The space $X$ is therefore the natural space for $z$,
corresponding to the energy space for the original map ${\bf u}$.
The facts
\begin{equation}
\label{eq:simple}
  z \in X \; \implies z \mbox{ continuous in } (0,\infty),
  \; z(0+) = z(\infty-) = 0,
  \; \mbox{and} \;
  \| z \|_{L^\infty} \lec \| z \|_X,
\end{equation}
follow easily from the change of variable $\rho^m = e^y$ and
Sobolev imbedding on $\R$ (see \cite{GKT}).

\subsection{Equation for the perturbation}
\label{sec:z}

The next step is to derive an equation for $z(\rho,t)$.
In terms of ${\bf v}(r,t)$, the Schr\"odinger map equation can be written as
\begin{equation}
\label{eq:v}
  {\bf v}_t = {\bf v} \times \left(
  {\bf v}_{rr} + \frac{1}{r} {\bf v}_r
  + \frac{m^2}{r^2} R^2 {\bf v} \right).
\end{equation}
Using~\eqref{eq:split}, we find
\begin{equation}
\label{eq:left}
  e^{-\al R} {\bf v}_t =
 [\dot \al R - s^{-1} \dot s \rho \pd_\rho](\bh+\bxi) + \bxi_t,
\end{equation}
\begin{equation}
\label{eq:right}
  s^2 e^{-\al R}({\bf v} \times M_r{\bf v}) =
 (\bh+\bxi) \times (M_\rho \bh + M_\rho\bxi),
\end{equation}
where
\[
  M_\rho := \pd_{\rho}^2 + \frac{1}{\rho} \pd_\rho
  + \frac {m^2}{\rho^2}R^2
\]
(and the right-hand sides are evaluated at $(\rho = r/s(t),t)$).

Consider first~\eqref{eq:right}.
Since $\Delta {\bf H} + |\nabla {\bf H}|^2 {\bf H}=0$ for
${\bf H} = e^{m \th R} \bh$, we have
\begin{equation}
\label{eq:Mh}
  M \bh = - 2\frac {m^2}{\rho^2} h_1^2 \bh,
\end{equation}
where $M = M_\rho$. Thus,
\begin{align*}
  \text{RHS of } \eqref{eq:right}
  &= \bh \times M \bxi +  \bxi \times (- 2\frac{m^2}{\rho^2} h_1^2 \bh)
  + \bxi \times M \bxi \\
  &= \bh \times (M +  2\frac {m^2}{\rho^2} h_1^2)\bxi + \bxi \times M \bxi.
\end{align*}
Keeping in mind~\eqref{eq:gamma2}, we write
\[
  \text{RHS of } \eqref{eq:right}
  =\bh \times (M + 2\frac {m^2}{\rho^2} h_1^2)(\bV^\rho(z)) + {\bf F_1}
\]
where ${\bf F_1} = \bh \times (M + 2\frac {m^2}{\rho^2} h_1^2)\gamma \bh + \bxi
\times M \bxi$ is the nonlinear part.
By~\eqref{eq:Mh}, we have
$M \gamma \bh = 2 \gamma_\rho \bh_\rho + (\cdots)\bh
= 2 \gamma_\rho \frac m\rho \hk + (\cdots)\bh$, and hence
\begin{equation}
\label{eq:F1}
  {\bf F_1} = -2 \gamma_\rho \frac m\rho h_1 \hj + \bxi \times M \bxi.
\end{equation}
Using $R^2 \hj = - \hj$,
$R^2 J^\bh \hj = h_1h_3 \bh -h_3^2 J^\bh \hj$,
$(J^\bh \hj)_\rho = -\frac{m}{\rho}h_1 \bh$, and
$(J^\bh \hj)_{\rho \rho} = -\frac{m^2}{\rho^2} h_1^2 J^\bh \hj -
(\frac{m}{\rho}h_1)_\rho \bh$ (all easy computations),
we find that the linear part can be rewritten as
\begin{align*}
  \bh \times (M +  2\frac {m^2}{\rho^2} h_1^2)(\bV^\rho(z)))
  &= - \bh \times [\bV^\rho(Nz)] = \bV^\rho(-iNz)
\end{align*}
where $N$ denotes the differential operator
$N := -\pd_\rho^2 - \frac 1\rho \pd_\rho + \frac{m^2}{\rho^2}(1-2h_1^2)$.

Because $\bxi_t = \bV^\rho(z_t) + \gamma_t \bh$,
\eqref{eq:left}--\eqref{eq:right} give
\[
  s^2 [\bV^\rho(z_t) + \gamma_t \bh]
  +[s^2 \dot \al R - s \dot s \rho \pd_\rho](\bh+\bxi)
  = \bV^\rho(-iNz) + {\bf F_1},
\]
or
\begin{equation}
\label{eq:z1}
  \bV^\rho(s^2 z_t + iN z) = {\bf F},
\end{equation}
where
\[
  {\bf F} := {\bf F_1} + [-s^2 \dot \al R +s \dot s \rho \pd_\rho ](\bh+\bxi)
  - s^2 \gamma_t \bh.
\]
Because the l.h.s. of~\eqref{eq:z1} is $\in T_\bh \bS^2$,
the r.h.s is also, and hence ${\bf F} \cdot \bh \equiv 0$.
We can re-write~\eqref{eq:z1} on the complex side
by applying $(V^\rho)^{-1}$:
\begin{equation}
\label{eq:z}
  i s^2 \frac{\p z}{\p t} = N z + i (\bV^\rho)^{-1}{\bf F},
  \quad\quad N = -\pd_\rho^2 - \frac 1\rho \pd_\rho +
  \frac{m^2}{\rho^2}(1-2h_1^2).
\end{equation}
This is the equation we sought for $z(\rho,t)$.

In order to see the form of the ``nonlinear'' terms
$(\bV^\rho)^{-1}({\bf F})$ more clearly, we compute
\[
  (\bV^\rho)^{-1}(R \bh(\rho)) = h_1(\rho), \quad
  (\bV^\rho)^{-1}(\rho \pd_\rho \bh(\rho)) = i mh_1,
\]
\[
  (\bV^\rho)^{-1}(P^{\bh(\rho)} R \bV^\rho(z)) = i z h_3, \quad
  (\bV^\rho)^{-1}(P^{\bh(\rho)} \rho\pd_\rho \bV^\rho(z))
  = \rho z_\rho,
\]
where $P^{\bh(\rho)}$ denotes the orthogonal vector projection onto
$T_{\bh(\rho)} \bS^2$.
Thus, using $\bh+\bxi = (1+\gamma)\bh + \bV^\rho(z)$,
\begin{equation}
\label{eq:nlterms}
  (\bV^\rho)^{-1}({\bf F})
  = [-s^2\dot \al + im s \dot s ] (1+\gamma) h_1
  - s^2\dot \al i z h_3 + s \dot s \rho z_\rho
  + (\bV^\rho)^{-1}(P^{\bh(\rho)} {\bf F_1}).
\end{equation}

\subsection{Orthogonality condition and parameter equations}
\label{sec:orth}

We have not yet specified $s(t)$ and $\alpha(t)$.
The main result of~\cite{GKT} says that if the energy is close to
$\Emin$, that is $\delta_1^2 := \sE({\bf u}) - \sE_{min} \ll 1$,
then there exist continuous $s(t) > 0$ and $\alpha(t) \in \R$ such that
$\| e^{m \theta R} \bxi \|_{\dot H^1} \lec \delta_1$ as long as $s(t)$
stays away from $0$. The choice of the parameters was simple and natural:
at each time $t$, $s(t)$ and $\alpha(t)$ were chosen so as to minimize
$\| e^{m \theta R} \bxi \|_{\dot H^1}$.
In this paper, we are forced into a different choice of $s(t)$
and $\alpha(t)$, as we shall now explain.

Supposing for a moment that $s(t) \equiv 1$, the linearized equation
for $z(\rho,t)$ can be read from~\eqref{eq:z}:
\begin{equation}
\label{eq:zlin}
  i \pd_t z = N z.
\end{equation}
The factorization
\begin{equation}
\label{eq:factor}
  N = L_0^* L_0, \quad\quad
  L_0 := \pd_\rho + \frac{m}{\rho}h_3 = h_1 \pd_\rho \frac{1}{h_1}
\end{equation}
(where the adjoint $L_0^*$ is taken in the $L^2(\rho d\rho)$
inner product)
shows that $ker N = span \{ h_1 \}$. In particular,~\eqref{eq:zlin}
admits the constant (in time) solution
$z(\rho,t) \equiv h_1(\rho)$.
Since we would like $z(\rho,t)$ to have some decay in time,
we must choose $s(t)$ and $\alpha(t)$ in such a way as to avoid
such constant solutions.
Since $N$ is self-adjoint in $L^2$, the natural choice is to work in the
subspace of functions $z$ satisfying
\begin{equation}
\label{eq:orth}
  (z, \; h_1)_{L^2}  = \int_0^\infty z(\rho) h_1(\rho)
  \rho d \rho \equiv 0,
\end{equation}
which is invariant under the linear flow~\eqref{eq:zlin}.

Recall, however, that the ``energy space'' for $z$ is
the space $X$ (defined in~\eqref{eq:X}).
Certainly the linear flow~\eqref{eq:zlin} does {\it not} preserve the
subspace $\bket{f \in X, \bka{f,h_1}_X =0}$
(since $N$ is not self-adjoint in $X$).
In fact, neither $z$ nor $h_1$ lies in $L^2$ in general.
The best we can do is
\[
  |(z, h_1)_{L^2}| =
  \left| \left( \frac{z}{\rho}, \rho h_1 \right)_{L^2} \right|
  \leq \| z \|_X \| \rho h_1 \|_{L^2}.
\]
So to make sense of~\eqref{eq:orth}, we require
\begin{equation}
\label{eq:restriction}
  \rho h_1(\rho) = \frac{2\rho}{\rho^m + \rho^{-m}} \in L^2(\rho d\rho),
\end{equation}
which only holds if $m \geq 3$.
This is one of the reasons we cannot handle
the small  $|m|$ cases in Theorem~\ref{thm:gkt2}.
The further restriction $m > 3$ is needed in
Proposition~\ref{prop:zbyq} to come.

In order to ensure condition~\eqref{eq:orth} holds for all
times $t$, it suffices to impose it initially,
and then ensure the time derivative of the inner-product
vanishes for all $t$.
Differentiating~\eqref{eq:orth} with respect to $t$, and using
Equations~\eqref{eq:z}, \eqref{eq:nlterms},
and~\eqref{eq:orth}, yields a system of ODEs for $s(t)$ and $\al(t)$:
\begin{equation}
\label{eq:ode}
  [s^2 \dot\al - im s \dot s ] (h_1,\, (1+\gamma)h_1)_{L^2}
  = (h_1, \ (\bV^\rho)^{-1}(P^{\bh(\rho)} {\bf F_1}) - s^2\dot \al i h_3 z
  + s \dot s \rho z_\rho)_{L^2}.
\end{equation}
The orthogonality condition~\eqref{eq:orth} is precisely
the one that ensures the terms linear in $z$ disappear
from~\ref{eq:ode}, and hence the
key property that $\dot s$ and $\dot \alpha$
be at least {\it quadratic} in $z$.
More precisely, the system~\eqref{eq:ode}
leads to the following estimate:
\begin{lemma}
\label{lem:odeest}
If $\| z \|_X \ll 1$, then
\[
  |s \dot{s}| + |s^2 \dot{\alpha}| \lec
  \norm{\frac{z}{\rho^2}}_{L^2}^2 + \norm{\frac{z_\rho}{\rho}}_{L^2}^2.
\]
\end{lemma}
\begin{proof}
Using
\[
\begin{split}
  |( h_1, h_3 z)| &\lec \| \rho h_1 \|_{L^2} \| z/\rho \|_{L^2}
  \lec \| z \|_X \ll 1 \\
  |(h_1, \rho z_{\rho})| & \lec \| \rho h_1 \|_{L^2} \| z_{\rho} \|_{L^2}
  \lec \| z \|_X \ll 1 \\
  |(h_1, \gamma h_1)| & \lec \| \rho^2 h_1^2 \|_{L^\infty} \|z/\rho\|_{L^2}^2
  \lec \| z \|_X^2 \ll 1,
\end{split}
\]
in~\eqref{eq:ode}, we arrive at
\begin{equation}
\label{eq:firstest}
  |s \dot{s}| + |s^2 \dot{\alpha}| \lec
  |(h_1, (\bV^\rho)^{-1}(P^{\bh(\rho)} {\bf F_1}))|.
\end{equation}
To finish the proof of the lemma, we will need to find
$(\bV^\rho)^{-1}(P^\bh {\bf F_1})$ explicitly.
Using the calculation of Lemma~\ref{Ap2:T2} in Appendix B,
we have
\[
\begin{split}
  (h_1,(\bV^\rho)^{-1} P^\bh {\bf F_1})_{L^2}
  = \int_0^\infty \bigg( & i (h_1)_\rho (- \gamma z_{\rho} + z \gamma_{\rho})
  + \frac{m}{\rho} h_1^2 (-2 \gamma_{\rho} - i z_2 (z_1)_{\rho} + i z_1 (z_2)_{\rho})
  \\
  & + \frac{m}{\rho}(h_1^2)_{\rho} (\gamma^2 - i z_2 z)
  + i \frac {m^2}{\rho^2}(2h_1^2-1) h_1 \gamma z \bigg) \rho d\rho.
\end{split}
\]
Now using the inequality \eqref{eq:simple}, together with
$(h_1)_{\rho} = -(m/\rho)h_1 h_3$, and the fact that $\rho^2
h_1(\rho)$ is bounded for $m \geq 2$, the estimate
\[
 |(h_1,(\bV^\rho)^{-1}(P^\bh {\bf F_1})_{L^2}|
 \lec \norm{\frac{z}{\rho^2}}_{L^2}^2 +
 \norm{\frac{z_\rho}{\rho}}_{L^2}^2
\]
follows. Together with~\eqref{eq:firstest}, this
completes the proof of Lemma~\ref{lem:odeest}.
\end{proof}

\donothing{
Since
\[
  \bxi = z_1 \hj + z_2 J^\bh \hj + \gamma \bh = (\xi_1,z_1,\xi_3)^T
\]
with $\xi_1 = -z_2 h_3 + \gamma h_1$ and
$\xi_3 = z_2 h_1 + \gamma h_3$, and
\[
  \hi = - h_3 J^\bh \hj + h_1 \bh, \quad \hk = h_1 J^\bh \hj + h_3 \bh,
\]
we have
\[
  R^2 \bxi = - z_1 \hj -\xi_1 \hi
  = -z_1 \hj + \xi_1 h_3 J^\bh \hj - \xi_1 h_1 \bh
\]
and
\[
  \bxi \times R^2 \bxi   = \begin{vmatrix}
  \hj & J^\bh \hj & \bh \\
  z_1 & z_2   & \gamma \\
  -z_1 & \xi_1 h_3 & - \xi_1 h_1
  \end{vmatrix}
  = - \xi_1 (z_2 h_1 + \gamma h_3)\hj + z_1 (\xi_1 h_1
  - \gamma)J^\bh \hj + z_1(\xi_1 h_3 + z_2) \bh.
\]
Thus
\[
  (V^\rho)^{-1}(P^{\bh(\rho)}(\bxi \times \frac {m^2}{\rho^2}R^2 \bxi)) =
  \frac {m^2}{\rho^2} [- \xi_1 (z_2 h_1 + \gamma h_3) + i z_1 (\xi_1 h_1
  - \gamma)].
\]
Because
\[
  \bh_\rho = \frac{m}{\rho} h_1 J^\bh \hj, \quad \Delta \bh =
  - \frac{m^2}{\rho^2}h_1^2 \bh - \frac{m^2}{\rho^2}h_1h_3 J^\bh \hj,
\]
we have
\begin{align*}
  \Delta (\gamma \bh) &= [(\Delta - \frac{m^2}{\rho^2}h_1^2) \gamma] \bh
  + (-\frac{m^2}{\rho^2}h_1 h_3\gamma + \frac{2m}{\rho}h_1 \gamma_\rho)J^\bh \hj, \\
  \Delta(z_2 \bh) \times \hj &= [(\Delta - \frac{m^2}{\rho^2}h_1^2) z_2]
  J^\bh \hj +(\frac{m^2}{\rho^2}h_1 h_3 z_2 - \frac{2m}{\rho}h_1 (z_2)_\rho) \bh.
\end{align*}
Hence, denoting $A := \Delta - \frac {m^2}{\rho^2}h_1^2$ and
$B := - \frac{m^2}{\rho^2}h_1 h_3 + \frac{2m}{\rho}h_1 \pd_\rho$,
\begin{align*}
  (V^\rho)^{-1} (P^{\bh(\rho)} (\bxi \times \Delta \bxi))  &=
   (V^\rho)^{-1}(P^{\bh(\rho)} \begin{vmatrix}
  \hj & J^\bh \hj & \bh \\
  z_1 & z_2   & \gamma \\
  \Delta z_1 &  A z_2+ B \gamma &  A \gamma - B z_2
  \end{vmatrix}) \\
  & = z_2(A \gamma - B z_2) - \gamma (A z_2 + B \gamma)
  + i \gamma \Delta z_1 - i z_1(A \gamma - B z_2)
\end{align*}
Finally, by \eqref{eq:F1}
\begin{align*}
  (V^\rho)^{-1} (P^{\bh(\rho)} {\bf F_1}) & = - 2 \frac{m}{\rho} h_1 \gamma_\rho +
  z_2(A \gamma - B z_2) - \gamma (A z_2+ B \gamma)
  + i \gamma \Delta z_1 - i z_1(A \gamma - B z_2) \\
  & \quad +  \frac{m^2}{\rho^2} [- \xi_1 (z_2 h_1 + \gamma h_3) +
  i z_1 (\xi_1 h_1 - \gamma)] \\
  & =-2 \frac{m}{\rho} h_1 \gamma_\rho + i z(Bz_2- \frac{m^2}{\rho^2}h_1 h_3 z_2)
  + O(z^3) \\
  & =-2 \frac m\rho h_1 \gamma_\rho + i z \frac{2m}{\rho} \pd_\rho (h_1 z_2)
  + O(z^3)
\end{align*}
Now we may estimate the inner product
$(h_1, (V^\rho)^{-1}(P^{\bh(\rho)} {\bf F_1}))$.
We have
\begin{itemize}
\item
$|(h_1, \frac{m}{\rho} h_1 \gamma_\rho)|
\lec \int \rho^2 h_1^2 |z/\rho^2| |z_\rho/\rho| \rho d\rho
\lec \| z/\rho^2 \|_{L^2} \|z_\rho/\rho\|_{L^2}$
\item
$|(h_1, \frac{z}{\rho} (h_1 z_2)_{\rho})|
\lec \int \rho^2 h_1^2 |z/\rho^2|(|z_\rho/\rho|+|z/\rho^2|)
\lec \| z/\rho^2 \|_{L^2}(\| z_\rho/\rho \|_{L^2} + \| z/\rho^2 \|_{L^2})$
\item cubic and quartic terms: after integration by parts in terms
involving two derivatives, using the simple inequality
\begin{equation}
\label{eq:simple}
  \| z \|_{L^\infty} \lec \| z \|_X
\end{equation}
from~\cite{GKT}, and the boundedness of $h_1(\rho)$ and
$\rho^2 h_1(\rho)$, these terms are all controlled by
\[
  \int h_1 |z| (|z_\rho|^2 + |z| |z_\rho|/\rho + |z|^2/\rho^2)
  \lec  \| z_\rho/\rho \|_{L^2}^2 + \| z/\rho^2 \|_{L^2}^2.
\]
\end{itemize}
Combining these estimates with~\eqref{eq:firstest}
completes the proof of Lemma~\ref{lem:odeest}.
}

\subsection{A nonlinear Schr\"odinger equation suited to estimates}
\label{sec:q}

We need to prove that $z(\rho,t)$ has some decay in time,
but the nonlinear Schr\"odinger-type equation~\eqref{eq:z}
is not suitable for obtaining such estimates, for at least two reasons.
Firstly, as remarked previously, the linearized equation has
constant solutions, and so the orthogonality condition~\eqref{eq:orth}
has to be explicitly used in order to get any decay whatsoever.
Secondly, and maybe more seriously, some of the nonlinear terms contain
derivatives (even two derivatives) of $z$, leading to a loss of regularity.
Fortunately, there is a neat way around these problems:
the {\it generalized Hasimoto transform} of \cite{CSU}
yields an equation without these difficulties,
as we now explain.

Let ${\bf u} = e^{m\theta R} {\bf v}(r) \in \Sigma_m$.
From~(\ref{eq:energy2}), it is clear that the tangent vector
\[
  {\bf W}(r) := {\bf v}_r(r) - \frac{m}{r} J^{\bf v} R {\bf v}(r)
  \;\; \in T_{{\bf v}(r)} \bS^2
\]
plays a distinguished role. In particular, ${\bf u}$ is a harmonic map
if and only if  ${\bf W} \equiv 0$.
Indeed, the Schr\"odinger map equation~(\ref{eq:sm}),
written in terms of ${\bf v}(r,t)$, can be factored as
\begin{equation}
\label{eq:factored}
  \frac{\pd {\bf v}}{\pd t} = J^{\bf v}
  [D^{\bf v}_r + \frac{1}{r} -\frac{m}{r} v_3] {\bf W}
\end{equation}
where
\[
  D^{\bf v}_r := P^{\bf v(r)} \p_r
\]
denotes the {\it covariant derivative} (with respect to $r$, along ${\bf v}$).
The idea is to write an equation for ${\bf W}$ in an appropriate intrinsic way.

Following~\cite{CSU}, let ${\bf e}(r) \in T_{{\bf v}(r)} \bS^2$
be a unit-length tangent field satisfying the ``gauge condition''
\begin{equation}
\label{eq:gauge}
  D^{\bf v}_r {\bf e} \equiv 0.
\end{equation}
Expressing ${\bf W}$ in the orthonormal frame
$\{ {\bf e}, \; J^{\bf v} {\bf e} \}$,
\[
  {\bf W} = q_1 {\bf e} + q_2 J^{\bf v} {\bf e},
\]
and using~\eqref{eq:factored}, and~\eqref{eq:gauge}, it is not
difficult to arrive at the following equation for the complex
function $q(r,t) := q_1(r,t) + iq_2(r,t)$:
\begin{equation}
\label{eq:q}
\begin{split}
  iq_t &= -(\p_r + \frac{m}{r} v_3)(\p_r + \frac{1}{r} - \frac{m}{r} v_3) q
  + S q \\
  &= (-\Delta_r + \frac{1}{r^2}((1 - m v_3)^2 + mr(v_3)_r))q
  + S q
\end{split}
\end{equation}
where the function $S(r,t)$ arises as $D^{\bf v}_t {\bf e} = S
J^{\bf v} {\bf e}$. From the curvature relation
\[
  [D_r, D_t] {\bf e} =
  -Re \left[ \left( \p_r + \frac{1}{r} - \frac{m}{r} v_3 \right) q \;
   \overline{ \left(q + \frac{m}{r} \nu \right) } \right] J^{\bf v}
  {\bf e},
\]
where $P^{\bf v(r)} \hk = \hk - v_3 {\bf v}
= \nu_1 {\bf e} + \nu_2 J^{\bf v} {\bf e}$,
we find
\begin{equation}
\label{eq:NL}
  S = Re \int_r^\infty
  \left( \p_\tau + \frac{1}{\tau} - \frac{m}{\tau} v_3(\tau,t) \right) q(\tau,t) \;
   \overline{ \left( q(\tau,t) + \frac{m}{\tau} \nu(\tau,t) \right) } d\tau.
\end{equation}
Thus the term in~\eqref{eq:q} involving $S$ is non-local and
nonlinear. We can simplify the expression for $S$ by integrating
by parts in the term involving $\partial_\tau q$, and using the
relation $\nu_r = -v_3(q + (m/r)\nu)$, to arrive at
\begin{equation}
\label{eq:NL2}
  S(r,t) = -\frac{1}{2} Q(r,t)
  + \int_r^\infty \frac{1}{\tau} Q(\tau,t) d\tau,
  \quad\;\;
  Q := |q|^2 + \frac{2m}{r} Re (\bar{\nu} q).
\end{equation}
Thus Equation~\eqref{eq:q} resembles a cubic nonlinear Schr\"odinger
equation, keeping in mind (a) there are non-local nonlinear terms,
and (b) it is not self-contained: the unknown map
${\bf v}(r,t)$ itself appears in several places
(including through $\nu$).
Furthermore, since
\[
  \delta_1^2 = \sE({\bf u}) - 4\pi m = \frac{1}{2}\| {\bf W} \|_{L^2}^2
  = \pi \| q \|_{L^2(rdr)}^2,
\]
we are dealing with a {\it small $L^2$-data problem} for
Equation~\eqref{eq:q} (even though the map ${\bf u}$
is not a small-energy map). This is what allows us the estimates we need.

Because of the fact (b) mentioned above, and in order to close the
estimate of Lemma~\ref{lem:odeest}, we need to
be able to control $z$ (and hence ${\bf v}$) in terms of $q$.
This is only possible if we have
a supplementary condition such as~\eqref{eq:orth}
(since $q=0$ just means ${\bf v}(r) = e^{\alpha R}h(r/s)$
for some $s$, $\alpha$). Parts of the proof of the following estimates
are a simple adaptation of the corresponding argument in~\cite{GKT},
where the orthogonality condition was somewhat different.
\begin{proposition}
\label{prop:zbyq} If $m \geq 3$ and~\eqref{eq:orth} holds, and if
$\| z \|_X \ll 1$, then for $2 \leq p < \infty$,
\begin{enumerate}
\item $ \norm{z_\rho}_{L^p} + \norm{\frac{z}{\rho}}_{L^p}
        \lec s^{1-2/p} \| q \|_{L^p}$
\item if $m > 3$, $\norm{\frac{z_{\rho}}{\rho}}_{L^2} +
        \norm{\frac{z}{\rho^2}}_{L^2}
        \lec s \norm{ \frac{q}{r} }_{L^2}.$
\end{enumerate}
\end{proposition}
\begin{proof}
The first observation is that, modulo nonlinear terms,
$q(r)$ is equivalent to $(1/s)(L_0 z)(r/s)$, where
$L_0 = \partial_\rho + \frac{m}{\rho} h_3(\rho)$.
Precisely,
\[
\begin{split}
  s {\bf W}(s \rho) =& \bV^{\rho}( L_0 z)
  + \frac{m}{\rho} z_1 (h_1 z_2 + h_3 \gamma) \hj \\
  & \quad +
  \frac{m}{\rho}(-h_1 z_1^2 + [h_3 z_2 - h_1(1+\gamma)]\gamma) J^{\bh}\hj
  + (\gamma_\rho + \frac{m}{\rho}[h_1 z_2 \gamma - h_3 |z|^2]) \bh.
\end{split}
\]
Using~\eqref{eq:simple},
it follows easily that for $2 \leq p \leq \infty$,
\[
\begin{split}
  \| L_0 z \|_{L^p} &\lec s^{1-2/p} \| q \|_{L^p}
  + (\| z \|_X + \| z \|_X^3)\| |z_\rho| + |z|/\rho \|_{L^p} \\
  \| \frac{1}{\rho} L_0 z \|_{L^2} &\lec s \| \frac{1}{r} q \|_{L^2}
  + (\| z \|_X + \| z \|_X^3)\| |z_\rho|/\rho + |z|/\rho^2 \|_{L^2}.
\end{split}
\]
In light of these estimates, and $\| z \|_X \ll 1$,
Proposition~\ref{prop:zbyq} follows from the following lemma.
\begin{lemma}
\label{lem:zbyq} For $m \geq 3$ and $z(\rho)$
satisfying~\eqref{eq:orth},
\begin{enumerate}
\item $\| z \|_X \lec \| L_0 z \|_{L^2}$ \item $\norm{ |z_\rho| +
\frac{|z|}{\rho} }_{L^p} \lec \| L_0 z \|_{L^p}$ for $2 \leq p <
\infty$ \item if $m > 3$, $\norm{ \frac{|z_\rho|}{\rho} +
\frac{|z|}{\rho^2} }_{L^2} \lec \norm{ \frac{L_0 z}{\rho}
}_{L^2}$.
\end{enumerate}
\end{lemma}
{\it Proof of the lemma.}
An estimate very similar to the first one here is proved in~\cite{GKT}
(only the orthogonality condition is different).
Here we prove the first and third statements together, by showing
\[
  \| |z_\rho|/\rho^b + |z|/\rho^{1+b} \|_{L^2}
  \lec \| L_0 z/\rho^b \|_{L^2}
\]
for $-1 \leq b \leq 1$.
If this is false, we have a sequence $\{ z_j \}$, with
\begin{equation}
\label{eq:norm1}
\begin{split}
  &\| (z_j)_\rho/\rho^b \|_{L^2}^2 + \| z_j/\rho^{1+b} \|_{L^2}^2 = 1,
  \\
  & \int z_j(\rho)h_1(\rho) \rho d\rho = 0, \\
  & \| L_0 z_j/\rho^b \|_{L^2} \to 0.
\end{split}
\end{equation}
It follows that, up to subsequence,
$z_j \to z^*$ weakly in $H^1$ and strongly in $L^2$
on compact subsets of $(0,\infty)$,
and that $L_0 z^* = 0$.
Hence $z^*(\rho) = C h_1(\rho)$ for some $C \in \C$.
Integration by parts gives
\[
  \| L_0 z_j/\rho^b \|_{L^2}^2 =
  \| (z_j)_\rho/\rho^b \|_{L^2}^2
  + m \int_0^\infty \frac{|z_j|^2}{\rho^{2b+2}}
    (m + 2bh_3(\rho) - 2mh_1^2(\rho)) \rho d\rho
\]
and so, defining $V(\rho) := m + 2bh_3(\rho) - 2mh_1^2(\rho)$,
we see that for any $\epsilon < 1/m$,
\[
  \limsup_{j \to \infty}
  m \int_0^\infty \frac{|z_j|^2}{\rho^{2b+2}} [V(\rho) - \epsilon] \rho d\rho
  \leq -m\epsilon.
\]
If $2|b| + \epsilon < m$
(which certainly holds under our assumptions $|b| \leq 1$ and $m > 3$),
then $\{ \rho \; | \; V(\rho) - \epsilon \leq 0 \}$
is a compact subset of $(0,\infty)$, and so
\[
  m \int_{V-\epsilon \leq 0} \frac{|C|^2 h_1^2(\rho)}{\rho^{2b+1}}
  [V(\rho) - \epsilon] \rho d\rho
  = \lim_{j \to \infty}
  m \int_{V - \epsilon \leq 0} \frac{|z_j|^2}{\rho^{2b+2}}
  [V(\rho) - \epsilon] \rho d\rho
  \leq -m \epsilon,
\]
which implies $C \not= 0$. Finally, for any $\epsilon' > 0$,
\[
\begin{split}
  0 &= \lim_{j \to \infty} \int_0^\infty z_j(\rho) h_1(\rho) \rho d \rho \\
    &= \int_{\epsilon'}^{1/\epsilon'} C h_1^2(\rho) \rho d \rho
    + \lim_{j \to \infty} \left( \int_0^{\epsilon'} + \int_{1/\epsilon'}^\infty \right)
    z_j(\rho) h_1(\rho) \rho d\rho.
\end{split}
\]
Since $\| z_j/\rho^{1+b} \|_{L^2} \leq 1$,
and $\rho^{1+b} h_1 \in L^2$
(this is precisely where we need $m > 3$, for $b=1$),
the last integrals are uniformly
small in $\epsilon'$, and we arrive at
\[
  0 = \int_0^\infty C h_1^2(\rho) \rho d \rho,
\]
contradicting $C \not= 0$.

%
%

We now prove the second statement.
First note that following the proof of Lemma 4.4 in~\cite{GKT},
the estimate
\begin{equation}
\label{eq:LpL2}
  \| |z_\rho| + |z|/\rho \|_{L^p} \lec
  \| L_0 z \|_{L^p} + \| L_0 z \|_{L^2}
\end{equation}
can be deduced from the $X$ estimate above (the case $b=0$).
Now fix a smooth cut-off function $\Phi(t)$ with
$\Phi(t)=1$ for $t \in [0,1]$, $\Phi(t)=0$ for $t \in
[2,\infty)$, and $\Phi_t(t)<0$ for $t \in (1,2)$.
Let $\phi(\rho) := \Phi(t)$ with $t = (\rho/s)^\beta$,
where $s \gg 1$ and $0 < \beta \ll 1$ are such that
\[
  \e_1 = \| \rho \phi_\rho(\rho) \|_{L^\infty}
  \lec \beta
\]
and
\[
  \e_2 := \norm{\rho [1 - \phi(\rho)] h_1(\rho)}_{L^{2}(\rho d\rho)}
  \le \norm{\rho h_1(\rho)}_{L^{2}((s,\infty),\rho d \rho)}
\]
are sufficiently small. Now using~\eqref{eq:orth},
\[
  |\int h_1 z \phi \rho d\rho| =
  |\int h_1 z (1-\phi) \rho d\rho|
  \le \e_2  \norm{z/\rho}_{L^{2}(\rho d\rho)}
\]
Observe that the proof of the $X$ estimate above
(and hence also of~\eqref{eq:LpL2}), works even if
$|\int h_1 z \rho d\rho| = o(1) \|z/\rho\|_{L^2}$,
and so provided $\e_2$ is sufficiently small,
we can apply~\eqref{eq:LpL2} to obtain
\begin{align*}
\norm{z/\rho}_p &\le \norm{z \phi /\rho}_p + \norm{z (1-\phi) /\rho}_p
\\
& \lec \norm{L_0(z \phi)}_p + \norm{L_0(z \phi)}_2
+\norm{z (1-\phi)/\rho}_p
\\
& \lec \norm{L_0(z \phi)}_p + \norm{z(1-\phi)/\rho}_p
\end{align*}
Now $1 - \phi$ is supported for $\rho \geq s \gg 1$,
and on this set $h_3(\rho) \geq 1/2$.
Then an easy adaptation of Lemma 4.2 in~\cite{GKT}
(using $m > 1$) yields
\[
  \norm{z(1-\phi)/\rho}_p \lec \norm{L_0(z(1-\phi))}_p,
\]
and hence
\begin{align*}
  \norm{z/\rho}_p &\lec
  \norm{L_0(z \phi)}_p + \norm{L_0(z(1-\phi))}_p \\
  & \lec \norm{L_0(z) \phi}_p + \norm{L_0(z)(1-\phi)}_p
  + \norm{z\phi_\rho}_p.
\end{align*}
Since $\norm{z \phi_\rho}_p \le \e_1 \norm{z/\rho}_{p}$, we conclude
\[
  \norm{z_\rho}_p + \norm{z/\rho}_p \le C
  \norm{L_0(z)}_p + C \e_1 \norm{z/\rho}_{p}.
\]
If $\e_1$ is small enough, the last term can be absorbed to the left side.

That completes the proof of the lemma, and hence of
Proposition~\ref{prop:zbyq}.
\end{proof}

Combining Proposition~\ref{prop:zbyq} with Lemma~\ref{lem:odeest}
leads to
\begin{corollary}
\label{cor:sdot} Under the conditions of
Proposition~\ref{prop:zbyq}, if $m > 3$,
\begin{equation}
\label{eq:sbyq}
  |s^{-1} \dot{s}| + |\dot \alpha| \lec \| q/r \|_{L^2}^2.
\end{equation}
\end{corollary}
This is our main estimate of the harmonic map parameters
$s(t)$ and $\alpha(t)$.

\subsection{Nonlinear estimates}

We can now use Proposition~\ref{prop:zbyq} to estimate the
nonlinear terms in~\eqref{eq:q}.
The idea is that from the splitting of Section~\ref{sec:splitting},
we expect $v_3(r,t) = h_3(r/s(t)) \; + \;$ ``small''.
We will ``freeze'' the scaling factor $s(t)$ at, say,
$s_0 := s(0)$ (and without loss of generality we will
rescale the solution so that $s_0 = 1$)
and treat the corresponding correction as a nonlinear term:
\begin{equation}
\label{eq:q2}
  iq_t + \Delta_r q - \frac{1 + m^2 - 2mh_3(r)}{r^2}q
  = Uq + Sq
\end{equation}
where
\[
  U := \frac{1}{r^2}
  [ m(v_3 - h_3)(m(v_3+h_3)-2) + mr((v_3)_r - (h_3)_r) ]
\]
(here we have used $r(h_3)_r = m h_1^2$ and $h_1^2 + h_3^2 = 1$),
and, recall from~\eqref{eq:NL2},
\[
  S(r,t) = -\frac{1}{2} Q(r,t)
  + \int_r^\infty \frac{1}{\tau} Q(\tau,t) d\tau,
  \quad\;\;
  Q := |q|^2 + \frac{2m}{r} Re (\bar{\nu} q).
\]
The next lemma estimates the r.h.s of~\eqref{eq:q2}
in various space-time norms.

\begin{lemma}
Provided~\eqref{eq:orth} holds, and $\| z \|_X \ll 1$, we have
\begin{equation}
\label{eq:nl1}
  \| r U q \|_{L^2_t L^2_x} \lec
  ((1 + \|s^{-1}\|_{L^\infty_t}) \| s - 1 \|_{L^\infty_t}
  + \| q \|_{L^\infty_t L^2_x} )
  \left\| \frac{q}{r} \right\|_{L^2_t L^2_x}
  + \| s^{-1} \|_{L^\infty_t}^{1/2} \| q \|_{L^4_t L^4_x}^2
\end{equation}
and
\begin{equation}
\label{eq:nl2}
  \| S q \|_{L^{4/3}_t L^{4/3}_x}
  \lec \| q \|_{L^4_t L^4_x}
  (\| q \|_{L^4_t L^4_x}^2 + \| q/r \|_{L^2_t L^2_x} ).
\end{equation}
\end{lemma}
\begin{proof}
Recall
\[
  v_3(r) = h_3(r/s) + \xi_3(r/s)
  = (1 + \gamma(r/s))h_3(r/s) + h_1(r/s)z_2(r/s),
\]
and set, as usual, $\rho = r/s$.
Estimate~\eqref{eq:nl1} follows from
$\| z \|_{L^\infty} \lec \| z \|_X$, the
estimates in Proposition~\ref{prop:zbyq}, and
\begin{itemize}
\item
$|h_3(r/s) - h_3(r)| = |\int_1^s \frac{d}{d \tau} h_3(r/\tau) d \tau|
= m |\int_1^s \frac{1}{\tau} h_1^2(r/\tau) d \tau|
\lec [\min(1,s)]^{-1} |s-1|$
\item
$r|[h_3(r/s)]_r - [h_3(r)]_r| = m|h_1^2(r/s)-h_1^2(r)|
\lec [\min(1,s)]^{-1}|s-1|$.
\end{itemize}
For estimate~\eqref{eq:nl2}, begin with
\[
  \| S q \|_{L^{4/3}_x L^{4/3}_t} \leq
  \| q \|_{L^4_t L^4_x} \| S \|_{L^2_t L^2_x}.
\]
Using the Hardy-type inequality
$\| \cdot \|_{L^2_x} \lec \| r \partial_r \cdot \|_{L^2_x}$ yields
\[
  \| S \|_{L^2_t L^2_x}
  \lec \| Q \|_{L^2_t L^2_x}
  \lec \| q \|_{L^4_t L^4_x}^2
  + \| \nu \|_{L^\infty_t L^\infty_x}
  \left\| \frac{q}{r} \right\|_{L^2_t L^2_x}.
\]
And since $|\nu| = |\hk - v_3 \bv| \lec 1$,
we arrive at~\eqref{eq:nl2}.
\end{proof}


\section{Dispersive estimates for critical-decay potentials
in two dimensions}
\label{sec:linear}

In order to establish any decay (dispersion) of solutions
of~\eqref{eq:q2}, we need good dispersive estimates for the linear part
\begin{equation}
\label{eq:linear}
  i q_t = -q_{rr} - \frac{1}{r}q_r + \frac{1}{r^2}(1+m^2-2mh_3) q
\end{equation}
This turns out to be a little tricky, since it is a ``borderline'' case
in two senses: the space dimension is two, and the potential
has $1/r^2$ behaviour both at the origin and at infinity, i.e.
\begin{equation}
\label{eq:Vas}
  \frac{1}{r^2}(1+m^2-2mh_3(r))
  \sim \left\{ \begin{array}{cc}
  \frac{(1+m)^2}{r^2} & \;\; r \to 0 \\
  \frac{(1-m)^2}{r^2} & \;\; r \to \infty
  \end{array} \right. .
\end{equation}

In this section we consider linear Schr\"odinger operators
like the one appearing on the r.h.s of~\eqref{eq:linear}.
More precisely, let
\begin{equation}
\label{eq:Vpos}
  H = -\Delta + \frac{1}{r^2} + V(r), \quad \;\;
  V \in C^\infty(0,\infty), \quad \;\;
  0 \leq r^2 V(r) \leq const.
\end{equation}
Such an operator is essentially self-adjoint on
$C_0^\infty(\R^2 \backslash \{0\})$, extends to a self-adjoint
operator on a domain $D(H)$ with
$C_0^\infty(\R^2 \backslash \{0\}) \subset D(H) \subset L^2(\R^2)$,
and generates a one-parameter unitary group
$e^{-i t H}$ such that for $\phi \in L^2$,
$\psi = e^{-i t H} \phi$ is the solution of the linear
Schr\"odinger equation $i \psi_t = H \psi$ with initial
data $\psi|_{t=0} = \phi$ (see, eg., \cite{RS2}).

Our goal is to obtain dispersive space-time ({\it Strichartz})
estimates for $e^{-itH}$ of the sort which hold for
the ``free'' ($H = -\Delta$) evolution:
\begin{equation}
\label{eq:free}
  \| e^{it\Delta} \phi \|_{L^r_t L^p_x}
  + \norm{ \int_0^t e^{i(t-s)\Delta} f(s) ds }_{L^r_t L^p_x}
  \lec \| \phi \|_{L^2} +
  \| f \|_{L^{\tilde{r}'}_t L^{\tilde{p}'}_x}
\end{equation}
where $(r,p)$ and $(\tilde{r},\tilde{p})$ are
{\it admissible} pairs of exponents:
\[
  (r,p) \mbox{ admissible } \quad \longleftrightarrow \quad
  1/r + 1/p = 1/2, \;\; 2 < r \leq \infty,
\]
and $p' = p/(p-1)$ denotes the H\"older dual exponent.
The {\it endpoint} case of \eqref{eq:free},
$(r,p) = (2,\infty)$, is known to be false in general,
but true for radial $\phi$ and $f$, save for the
``double endpoint'' case $r = \tilde{r} = 2$ (\cite{Tao}).

Perturbative arguments to extend estimates like~\eqref{eq:free} to
Schr\"odinger operators with potentials (in general one
has to include a projection onto the continuous spectral
subspace in order to avoid bound states, which do not disperse)
cannot work for borderline behaviour like~\eqref{eq:Vas}. 
Fortunately, the problem of obtaining dispersive estimates
when the potential has this critical fall-off (and singularity)
is taken up in a recent series of papers by
Burq, Planchon, Stalker, and Tahvildar-Zadeh
(see in particular~\cite{BPST1,BPST4}).
In place of a perturbative argument, the authors
make a repulsivity assumption on the potential
(which, in particular, rule out bound states), and
prove more-or-less directly -- by identities --
that solutions have some time decay,
in a spatially-weighted space-time sense
(a {\it Kato smoothing} - type estimate).
This approach is ideally suited to our present problem:
the operator appearing in~\eqref{eq:linear}
satisfies the following repulsivity property: when
written in the form~\eqref{eq:Vpos},
\begin{equation}
\label{eq:Vrep}
  -r^2 (r V(r))_r + 1 \geq \nu \quad
  \mbox{ for some } \;\; \nu > 0.
\end{equation}

We cannot rely directly on the results of ~\cite{BPST1,BPST4} here.
The paper~\cite{BPST1} considers only potentials $(const)/r^2$,
while the results of~\cite{BPST4} hold in dimension $\geq 3$ only,
and do not immediately extend to dimension two for two reasons:
one is the failure of the Hardy inequality, and the other is
the failure of the double-endpoint Strichartz estimate
(even for radial functions).
However, we can recover the argument from~\cite{BPST4} by
exploiting the radial symmetry of our functions to avoid the
Hardy inequality, and we can avoid the use of the double-endpoint
Strichartz estimate by following the approach of~\cite{BPST1},
which in turn follows~\cite{RodSchl}.

\begin{theorem}
\label{thm:dispersive}
Suppose the Schr\"odinger operator $H$ satisfies the
conditions~\eqref{eq:Vpos} and~\eqref{eq:Vrep}.
Let $\phi = \phi(r)$ be radially symmetric.
Then for any admissible pair $(r,p)$, we have
\begin{equation}
\label{eq:hom}
  \| e^{-itH} \phi \|_{L^r_t L^p_x} +
  \left\| \frac{1}{|x|} e^{-itH} \phi \right\|_{L^2_t L^2_x}
  \lec \| \phi \|_{L^2}.
\end{equation}
If $f = f(r,t)$ is radially symmetric,
and $(\tilde{r}, \tilde{p})$ is another admissible pair, then
\begin{equation}
\label{eq:inhom1}
  \left\| \int_0^t e^{-i(t-s)H} f(x,s) ds \right\|_{L^r_t L^p_x} +
  \left\| \frac{1}{|x|} \int_0^t e^{-i(t-s)H} f(x,s) ds \right\|_ {L^2_t L^2_x}
  \lec  \min \bke{ \| f \|_{L^{\tilde{r}'}_t L^{\tilde{p}'}_x},
  \| |x| f \|_{L^2_t L^2_x} }.
\end{equation}
\end{theorem}

\begin{remark}
In \cite{BPST4}, the single endpoint Strichartz estimate
(\eqref{eq:hom} with $r=2$) is also obtained for dimensions $\geq 3$.
In two dimensions, though it holds in the {\it free}, radial
case, we do not know if it holds for our operators.
However, it is {\it essential} to the present
paper to have an estimate with $L^2_t$ decay
($L^r_t$ with $r > 2$ is simply not enough -- see the next section).
Our way around this problem is to use the above {\it weighted}
$L^2_t L^2_x$ estimate that arises naturally in the
approach of~\cite{BPST4}.
\end{remark}

\begin{proof}
Parts of the proof are perturbative, so we
identify a reference operator:
\[
  H = -\Delta + \frac{1}{r^2} + V =: H_0 + V.
\]
Note that $H_0 = -\Delta + \frac{1}{r^2}$ satisfies the
'usual' Strichartz estimates
(those satisfied by $-\Delta$ as in~\eqref{eq:free} above)
on radial functions,
since $H_0$ is simply $-\Delta$ conjugated by $e^{i \theta}$
when acting on such functions.

{\bf Step 1.}
Following~\cite{BPST4}, we begin with weighted resolvent estimates.

\begin{lemma}
\label{lem:resolvent}
For $f=f(r)$ radial,
\begin{equation}
\label{eq:res}
  \sup_{\mu \not\in \R}
  \left\| \frac{1}{|x|} (H - \mu)^{-1} f \right\|_{L^2(\R^2)}
  \lec \| |x| f \|_{L^2(\R^2)}.
\end{equation}
\end{lemma}

\noindent
{\it Proof of Lemma}.
We can assume $f \in C_0^\infty(0,\infty)$, with the
lemma then following from a standard density argument.
Set $u := (H - \mu)^{-1} f$ so that
$(H - \mu) u = f$,
and note that $u = u(r)$ is radial, since $f$ is.
To avoid the use of the Hardy inequality in~\cite{BPST4},
we change variables from $u(r)$ to
\[
  v(x) := e^{i \theta} u(r)
\]
and use $|\nabla v|^2 = |u_r|^2 + \frac{1}{r^2} |u|^2$, so
\begin{equation}
\label{eq:newhardy}
  \left\| \frac{v}{|x|} \right\|_{L^2} \lec
  \| v \|_{H^1}
\end{equation}
In terms of $v$, the equation for $u$ becomes
\begin{equation}
\label{eq:resv}
  (-\Delta + V - \mu) v = \tilde{f}
\end{equation}
where $\tilde{f}(x) := e^{i \theta} f(r) \in L^2$,
and so $v \in D(-\Delta + V) \subset H^2$.
The proof of Lemma~\ref{lem:resolvent}
now follows precisely the corresponding proof in~\cite{BPST4},
using $-d^2/d \theta^2 \geq 1$ on functions of our
form $e^{i \theta} f(r)$, and
with~\eqref{eq:newhardy} (rather than Hardy) providing
$v/|x| \in L^2$ where needed. $\Box$

\smallskip

{\bf Step 2.}
As in \cite{BPST4}, the next step is to invoke~\cite{Kato} to conclude that
the resolvent estimate~\eqref{eq:res} implies the
following ``Kato smoothing'' weighted-$L^2$ estimate for the propagator:
for $\phi = \phi(r)$,
\begin{equation}
\label{eq:hom1}
  \left\| \frac{1}{|x|} e^{-itH} \phi \right\|_{L^2_t L^2_x}
  \lec \| \phi \|_{L^2}.
\end{equation}
This is one part of~\eqref{eq:hom}.
Note that the reference operator $H_0$ also satisfies the
weighted estimate~(\ref{eq:hom1})
(a fact which follows from the same argument).
Another direct consequence of the resolvent estimate~\eqref{eq:res}
is the inhomogeneous version of~\eqref{eq:hom1},
\begin{equation}
\label{eq:inhomx}
  \left\| \frac{1}{|x|} \int_0^t e^{-i(t-s)H} f(\cdot,s) ds \right\|_ {L^2_t L^2_x}
  \lec  \| |x| f \|_{L^2_t L^2_x},
\end{equation}
which is one part of~\eqref{eq:inhom1}.
The estimate~\eqref{eq:inhom1} is probably standard, but we
did not see a proof, and so supply one in Section~\ref{sec:weighted}.
\smallskip

{\bf Step 3.}
Next we establish more of the inhomogeneous
estimates in~\eqref{eq:inhom1},
but first for the reference operator $H_0$.
Since we do not have the double-endpoint Strichartz
estimate available, we now depart from~\cite{BPST4}
and henceforth follow~\cite{BPST1} (which in turn relies partly
on~\cite{RodSchl}).
Note that by~\eqref{eq:hom1} for $H_0$, for any $\psi \in L^2_x$,
\[
\begin{split}
  ( \psi, \int_0^\infty e^{isH_0} f(\cdot, s) ds )_{L^2_x}
  &= \int_0^\infty ds ( e^{-isH_0} \psi, f(\cdot, s) )_{L^2_x} \\
  & \leq \norm{\frac{1}{|x|} e^{-isH_0} \psi}_{L^2_t L^2_x}
  \| |x| f \|_{L^2_t L^2_x}
  \lec \|\psi\|_{L^2} \| |x| f \|_{L^2_t L^2_x},
\end{split}
\]
yielding
\[
  \norm{\int_0^\infty e^{isH_0} f(\cdot, s) ds}_{L^2_x}
  \lec  \| |x| f \|_{L^2_t L^2_x},
\]
and hence by the Strichartz estimates for $H_0$,
for $(r,p)$ admissible,
\[
\begin{split}
  \norm{\int_0^\infty e^{-i(t-s)H_0} f(\cdot, s) ds}_{L^r_t L^p_x}
  &= \norm{ e^{-itH_0} \int_0^\infty e^{isH_0} f(\cdot, s)
  ds }_{L^r_t L^p_x} \\
  & \lec \norm{ \int_0^\infty e^{isH_0} f(\cdot, s) ds }_{L^2_x}
  \lec  \| |x| f \|_{L^2_t L^2_x}.
\end{split}
\]
Finally, the required estimate
\begin{equation}
\label{eq:inhom0}
  \norm{ \int_0^t e^{-i(t-s)H_0} f(\cdot, s) ds }_{L^r_t L^p_x}
  \lec \| |x| f \|_{L^2_t L^2_x}
\end{equation}
follows from a general argument of Christ-Kiselev
(\cite{CK}, and see~\cite{BPST1}).

\smallskip

{\bf Step 4.}
To obtain the remaining part of~\eqref{eq:hom}
(the Strichartz estimate), we use~\eqref{eq:hom1},
and~\eqref{eq:inhom0}, in a perturbative argument.
We have
\[
  e^{-itH} \phi = e^{-itH_0} \phi
  + i\int_0^t e^{-i(t-s)H_0} V e^{-isH} \phi ds,
\]
and so for $(r,p)$ admissible,
\[
\begin{split}
  \| e^{-itH} \phi \|_{L^r_t L^p_x} & \lec
  \| \phi \|_{L^2} + \| |x| V e^{-isH} \phi \|_{L^2_t L^2_x} \\
  & \leq \| \phi \|_{L^2} + \| |x|^2 V \|_{L^\infty}
    \norm{ \frac{1}{|x|} e^{-isH} \phi }_{L^2_t L^2_x} \\
  &\lec \| \phi \|_{L^2}.
\end{split}
\]
This finishes the proof of~\eqref{eq:hom}.

\smallskip

{\bf Step 5.}
It remains to prove the rest of the inhomogeneous
estimates in~\eqref{eq:inhom1}.
But given~\eqref{eq:hom}, these follow again from the
argument used in Step 3.

That completes the proof of Theorem~\ref{thm:dispersive}.
\end{proof}

\begin{corollary}
\label{cor:Str}
If $m \geq 2$, the estimates~\eqref{eq:hom} and~\eqref{eq:inhom1}
hold for the operator
\[
  H := -\Delta + \frac{1}{r^2}(1+m^2-2mh_3)
\]
coming from the Schr\"odinger map problem.
\end{corollary}

\begin{proof}
We have
\[
  \frac{1}{r^2}(1+m^2-2mh_3) = \frac{1}{r^2} + V(r);
  \quad\quad
  V(r) = \frac{m}{r^2}(m-2h_3(r)).
\]
So for $m \geq 2$,
\[
  (m+1)^2 \geq 1 + r^2 V(r) \geq (m-1)^2 \geq 1,
\]
and
\[
  1 - r^2(rV)_r = 1 + m(m-2h_3(r)+2mh_1^2(r))
  \geq 1 + m(m-2) \geq 1.
\]
Thus the conditions~\eqref{eq:Vpos} and~\eqref{eq:Vrep}
both hold with $\nu = 1$.
\end{proof}


\section{Proof of the main theorem}
\label{sec:proof}

Let $\bu \in C([0,T_{max}); \Sigma_m)$ be the 
solution of the Schr\"odinger map equation~\eqref{eq:sm} with
initial data $\bu_0$ (given by Theorem~\ref{lwp:T1}).
Energy is conserved:
\[
  \sE(\bu(t)) = \sE(\bu_0) = 4\pi m + \delta_1^2.
\]
We begin by splitting the initial data $\bu(0)$,
using the following lemma, which is proved in Section~\ref{sec:split}:
\begin{lemma}
\label{lem:split}
If $m \geq 3$,and if $\delta$ is sufficiently small, then for any
map $\bu \in \Sigma_m$ with $\sE(\bu) \leq 4\pi m + \delta^2$,
there exist $s > 0$, $\al \in \R$,
and a complex function $z(\rho)$ such that
\begin{equation}
\label{eq:uform}
  \bu(r,\theta) = e^{[m\theta + \alpha]R}
  [ (1 + \gamma(r/s))\bh(r/s) + \bV^{r/s}(z(r/s)) ]
\end{equation}
with $z$ satisfying~\eqref{eq:orth}; i.e.,
\begin{equation}
\label{eq:orth2}
  \int_0^\infty z(\rho) h_1(\rho) \rho d \rho = 0,
\end{equation}
and $\| z \|_X^2 \lec \sE(\bu) - 4\pi m$.
\end{lemma}

Invoking the lemma, we have
\[
  \bu_0 = e^{[m\theta + \alpha_0]R}
  [ (1 + \gamma_0(r/s_0))\bh(r/s_0) + \bV^{r/s_0}(z_0(r/s_0)) ]
\]
with $z_0$ satisfying the orthogonality
condition~\eqref{eq:orth}, and
\[
  \| z_0 \|_X \lec \delta_1 \ll 1.
\]
Now rescale, setting
\[
  \hat{\bu}(x,t) := \bu(s_0 x, s_0^2 t).
\]
Then $\hat{\bu}$ is another solution of the Schr\"odinger
map equation~\eqref{eq:sm}, and
\[
  \hat{\bu}(x,0) = e^{[m\theta + \alpha_0]R}
  [ (1 + \gamma_0(r)) \bh(r) + \bV^{r}(z_0(r)) ].
\]
Let $q(r,t)$ be the complex function derived
from the Schr\"odinger map $\hat{\bu}$, as in
Section~\ref{sec:q}.

Suppose $(r,p)$ is an admissible pair of exponents.
Define a spacetime norm $Y$ by
\[
  \| q \|_Y :=
  \| q \|_{L^\infty_t L^2_x \cap L^4_t L^4_x \cap L^r_t L^p_x} +
  \left\| \frac{q}{r} \right\|_{L^2_t L^2_x}.
\]
As long as $\| z \|_X \lec \| q \|_{L^2_x}$ remains
sufficiently small, Corollary~\ref{cor:Str} together with
estimates \eqref{eq:nl1}-~\eqref{eq:nl2} yields
\begin{equation}
\label{eq:qestimates}
  \| q \|_Y \lec
  \| q(0) \|_{L^2} +
  \Big[ (1 + \|s^{-1}\|_{L^\infty_t})\|s-1\|_{L^\infty_t}
  + (1 + \|s^{-1}\|_{L^\infty_t})\| q \|_Y
  + \| q \|_Y^2 \Big] \|q \|_Y .
\end{equation}
We also have
\[
  \hat{\bu} = e^{[m\theta + \alpha(t)]R}
  [ (1 + \gamma(r/s(t),t))\bh(r/s(t)) + \bV^{r/s(t)}(z(r/s(t),t)) ],
\]
with $z(\rho,t)$ satisfying~\eqref{eq:orth},
$s(0) = 1$, $\alpha(0) = \alpha_0$, and, by Corollary~\ref{cor:sdot},
$s(t) \in C([0,T);\R^+)$ and $\alpha(t) \in C([0,T);\R)$, with
\begin{equation}
\label{eq:sestimates}
  \|s^{-1} \dot s\|_{L^1_t}
  + \| \dot \alpha \|_{L^1_t}
  \lec \| q \|_Y^2.
\end{equation}
Taking $\|q(0)\|_{L^2} \lec \delta_1$ sufficiently small,
the estimates~\eqref{eq:qestimates} and~\eqref{eq:sestimates}
yield
\begin{equation}
\label{eq:closed}
  \| q \|_Y \lec \delta_1, \quad \quad
   \|s^{-1} \dot s\|_{L^1_t} + \|\dot \alpha\|_{L^1_t} \lec \delta_1^2
\end{equation}
(and in particular, $\| z \|_X \ll 1$ continues to hold).
Since
\[
  | \nabla[\hat{\bu} -  e^{[m\theta + \alpha(t)]R} \bh(r/s(t))] |
  \lec \frac{1}{s}( |z_\rho| + |z/\rho| )( 1 + |z| ),
\]
the estimates of Proposition~\ref{prop:zbyq} give
\begin{equation}
\label{eq:spacetime}
  \| \nabla[\hat{\bu} -  e^{[m\theta + \alpha(t)]R} \bh(r/s(t))] \|
  _Y \lec \| q \|_Y \lec \delta_1.
\end{equation}
Estimate~\eqref{eq:closed} shows:
(a) that $s(t) \geq const > 0$, and hence, by
Corollary~\ref{lwp:T1b}, we must have $T_{max} = \infty$;
(b) that
\[
  s(t) \to s_\infty \in (1-c\delta_1^2,1+c\delta_1^2), \quad\quad
  \alpha(t) \to \alpha_\infty \in (\alpha_0-c\delta_1^2,\alpha_0+c\delta_1^2)
\]
as $t \to \infty$.

Finally, undoing the rescaling,
$\bu(r,t) = \hat{\bu}(r/s_0,t/s_0^2)$,
yields the estimates of Theorem~\ref{thm:gkt2}.
$\Box$


\donothing{
\section{Appendix: a lemma}

Let $\frac 1r L^2_{rad}$ denote the space of all functions $ f(r)$
with $rf(r)\in L^2_{rad}$.

Recall $N_0 = - \Delta _r + \frac {m^2}{r^2}$.

\bigskip

\begin{lemma} Let $m \ge 1$ be any positive integer.

(i) For all $ f \in \frac 1r L^2_{rad}$, there exists $g \in X$ so
that $N_0 g = f$.

(ii) Consider $N_0$ as an operator on $X$. Then $X \cap \frac 1r L^2_{rad}
\subset Ran(N_0)$.

\end{lemma}

\bigskip

Proof.
It suffices to show (i) since it implies (ii).
Consider real-valued functions on $(0,\infty)$.  The general solution
of the ODE $N_0g = 0$ is $g(r)=a r^m + b r^{-m}$.  The ODE $N_0g =f $
for any function $f(r)$ can be solved by the method of variations of
parameters:
\[
g(r) = r^m a(r) + r^{-m} b(r), \quad
g'(r) = \frac mr [r^m a(r) - r^{-m} b(r)],
\]
where
\[
a(r)=- \frac 1{2m}\int_{r_1}^r s^{1-m} f(s) ds , \quad
b(r)=\frac 1{2m}\int_{r_2}^r s^{1+m} f(s) ds ,
\]
for some constants $r_1,r_2$. We will choose $r_1 = \infty$ and $r_2 = 0$.
The condition $g \in X$, i.e., $g',g/r \in L^2_{rad}$, is equivalent to
\begin{equation}\label{eq1}
r^{m-1}\int_{\infty}^r s^{1-m} f(s) ds\in L^2_{rad} \quad \text{and} \quad
r^{-m-1}\int_{0}^r s^{1+m} f(s) ds  \in L^2_{rad}.
\end{equation}
Let $F(r) =\int_{\infty}^r s^{1-m} f(s) ds$, which vanishes as $ r \to
\infty$. (We may assume $f \in C^1_c(0,\infty)$ by density.) The first
of \eqref{eq1} is true if
\[
|r^{m-1}F|_{L^2_{rad}} \lesssim  |r^{m}F'|_{L^2_{rad}} =| rf|_{L^2_{rad}} ,
\]
but this follows from Hardy's inequality in $\R^{2m+2}$.  The second
of \eqref{eq1} follows from the first by duality: For any $k(r) \in
L^2_{rad}$,
\begin{align*}
&\int_0^\infty k(r) r^{-m-1}\int_{0}^r s^{1+m} f(s) ds \,r dr
=\int_0^\infty  sf(s) s^{m-1} \int_{s}^\infty r^{1-m} \frac{k(r)}r  dr\, s ds
\\
& \lesssim |sf(s)|_{L^2_{rad}} \cdot  |k(r)|_{L^2_{rad}}.
\end{align*}
}


\section{Appendix: local wellposedness}
\label{sec:lwp}

In this appendix we prove Theorem \ref{lwp:T1} and Corollary
\ref{lwp:T1b} on the local wellposedness of the Schr\"odinger flow
\eqref{eq:sm} when the data $\bu_0 \in \Sigma_m$ has energy
$\sE(\bu_0) = 4 \pi m + \de_0^2$
close to the harmonic map energy, $0<\de_0 \le \de
\ll 1$.  In subsection \ref{lwp:S1} we show that $z$
(and hence $\bu$)
can be reconstructed from $q$, $s$, and $\al$.  This subsection is
time-independent.  In subsection \ref{lwp:S2} we set up the equations
for the existence proof. In subsection \ref{lwp:S3} we show that we
have a contraction mapping, and complete the proof of Theorem
\ref{lwp:T1} and Corollary \ref{lwp:T1b}.  In subsection \ref{lwp:S4}
we discuss the small energy case.

Recall the decomposition $ \bu(r,\th) =e^{m \th R} \bv(r)$ and
\begin{equation}\label{lwp-05}
\bv(r) =  e^{ \al R}\, [\bh(\rho) + \bxi(\rho)]
= e^{ \al R}
\mat{(1+\gamma)h_1 - h_3 z_2 \\ z_1 \\ (1+\gamma)h_3 +h_1 z_2}(\rho),
\end{equation}
where $\rho = r/s$, $\bxi = z_1 \hj + z_2 h\times \hj + \gamma bh$ and
$\gamma = \sqrt{ 1- |z|^2} -1$.  The time-dependence of
$\bu,\bv,\bxi,\al,s$ and $\gamma$ has been dropped from \eqref{lwp-05}. The
equation $D_r \be =0$ is equivalent to
\begin{equation}
\be _r=-(\bv_r \cdot \be )\bv.
\label{lwp-9}
\end{equation}
Recall $q\be = \bv_r - \frac{m}{r}J^{\bv}R\bv$ with
$\nu \be = J^{\bv} R \bv = \hk - v_3 \bv$.
By substituting in \eqref{lwp-05} and using $L_0 \bh = \frac mr
\hk$, $q \be$ should satisfy
\begin{equation} \label{lwp-10}
s e^{-\al R}q \be(r)
= (L_0 z)(\rho) \hj + {\bf G_0}(z)(\rho), \quad \rho=r/s,
\end{equation}
where
\begin{equation}\label{lwp-15}
{\bf G_0}(z)(\rho) := s e^{-\al R}[ \bv_r - \frac{m}{r}(\hk - v_3 \bv)] -(L_0 z)\hj
=\gamma_\rho \bh + \frac m\rho (\gamma \hk + \gamma h_3 \bh + \xi_3 \bxi)
\end{equation}
and $\norm{{\bf G_0}(z)}_{L^2} \lec \nrm{z}_X^2$ when $\nrm{z}_X \ll 1$.
In other words, $q$ is rescaled $L_0 z$, plus error.

In this Appendix, we will choose a different orthogonality condition
for $z$, instead of \eqref{eq:orth}. Specifically, we choose the
unique $s$ and $\al$ so that
\begin{equation} \label{lwp:orth}
\bka{h_1, z}_X = 0.
\end{equation}
(Recall $\bka{f,g}_X = \int_0^\infty (\bar f_r g_r +
\frac{m^2}{r^2}\bar fg )rdr$.)
The condition \eqref{lwp:orth} makes sense for all $m \not = 0$ and
suffices for the proof of local wellposedness. In contrast,
\eqref{eq:orth} makes sense only if $|m| \ge 3$, but is
necessary for the study of the time-asymptotic behavior.
In \cite[Sect.~2]{GKT}, we
chose $s$ and $\al$ to minimize $\norm{\bu - e^{(m \th +\al)R}
\bh(\cdot/s) }_{\dot H^1}$.  The resulting equations in
\cite[Lem.~2.6]{GKT} are $ \bka{h_1, z_1}_X = 0$ and $\bka{h_1, z_2}_X
= \int_0^\infty \frac {4m^2}{\rho^2} h_1^2h_3 \gamma(\rho)\rho
d\rho$. The condition \eqref{lwp:orth} is similar but has no error
term. The unique choice of $s$ and $\al$ can be proved by implicit
function theorem, similar to the proof for Lemma~\ref{lem:split},
and is skipped. It is important to point out, however,
that the parameter $s$ used here, though not the 
same as $s(\bu)$ defined in~\eqref{eq:sdef}-\eqref{eq:sdef2},
is nonetheless comparable: $s = s(\bu)(1 + O(\delta_0^2))$
(this comes immediately from the implicit function theorem
argument). Thus we can state the local well-posedness result
(Theorem~\ref{lwp:T1}) in terms of $s(\bu_0)$.


\subsection{Reconstruction of $z$ and $\bu$ from $q$, $s$, and $\al$}
\label{lwp:S1}

In this subsection all maps are time-independent.  For a given map $\bu
= e^{m \th R}\bv(r) \in \Sigma_m$ with energy close to
$4 \pi m$, we can define $s,\al,z$ and $q$. The
three quantities $s,\al,z$ determine $\bu$, and hence $q$.
Conversely, as will be done in Lemma \ref{lwp:T3} of this subsection, we can recover
$z$ and $\bu$ if $s,\al$ and $q$ are given, assuming that
$\norm{q}_{L^2} \le \de$.  Before that we first prove difference
estimates for $\de \be$ in Lemma \ref{lwp:T2}.

For given $s>0$, $\al \in \R$ and $z \in X$ small, we define $\bv(r) =
\BV(z,s,\al)(r)$ by \eqref{lwp-05}, and
$\be(r) = \BE(z,s,\al)(r)$ by the ODE
\begin{equation} \label{lwp110}
\be(z)(0) = e^{\alpha R} \hj, \quad \be_r = -(\bv_r \cdot \be
)\bv, \quad \text{where } \bv = \BV(z,s,\al) .
\end{equation}
Also denote 
$\BE(z) = \BE(z,1,0)$. Simple comparison shows
\begin{equation} \label{lwp115}
\BE(z,s,\al)  = e^{\al R}\BE(z^s), \quad z^s(r) := z(r/s).
\end{equation}

\begin{lemma} \label{lwp:T2}
Suppose $z_l \in X$, $l=a,b$, are given with $\norm{z_l}_X$
sufficiently small. Let $\de z := z_a -z_b$,
$\de \bv := \BV(z_a,1,0)
-\BV(z_b,1,0)$, and $\de \be :=\BE(z_a) -\BE(z_b) $. Then
\[
  \norm{\de \bv}_{X} + \norm{\de \be}_{L^\infty}
  \lec \norm{\de z}_{X}.
\]
\end{lemma}

\begin{proof}
Note
\begin{equation}\label{lwp120}
\norm{\bh_r}_{L^2(rdr)}\leq C; \quad
\| \bxi^l \|_{X} \lec \| z^l \|_{X} + \| z^l \|_{X}^2, \quad
l=a,b.
\end{equation}
Since $\de \bv = \de \bxi = (\de z)\hj + (\de \gamma) \bh$,
\begin{equation}\label{lwp125}
\| \de \bv \|_X + \| \de \bxi \|_X
\lec  (1 + \| z_a\|_X + \| z_b\|_X  ) \|\de z \|_X
\lec   \|\de z \|_X.
\end{equation}

For $\de \be$, write $\de \be = (\de e_1, \de e_2,\de e_3)$ and
\begin{equation}\label{lwp130}
\de e_{j,r} = - (\de \bxi_r \cdot e_a)v_{a,j}
- (\bv_{b,r} \cdot \de \be)v_{a,j} - (\bv_{b,r}
\cdot \be_b)\de \bxi_j, \quad j=1,2,3.
\end{equation}

First consider $\de e_2$. Integrate in $r$. Using \eqref{lwp120},
\eqref{lwp125}, $v_{a,2} = z_{a,1}$, and $v_{l,r} \in L^2(rdr)$,
\begin{align}
\nonumber
|\de e_2(\tau)| &\lec \int_0^\tau \bke{|\de \bxi_r \cdot \be_a)\frac{z_a}r|
+ | (\bv_{b,r} \cdot \de \be)\frac {z_a}r|
+ | (\bv_{b,r} \cdot \be_b)\frac{\de z_1}r|} \,r d r
\\
&\lec (1+\max_{l=a,b} \norm{z_l}_X) \norm{\delta z}_{X}
+\max_{l=a,b} \norm{z_l}_X \norm{\delta \be }_{L^{\infty}}.
\end{align}

Next we consider $\delta e_1$ and $\delta e_3$.  Equations
\eqref{lwp130} for $j=1,3$ can be written as a vector equation for
$x=( \de e_1 , \de e_3)^T$:
\begin{equation} \label{lwp140}
x_r =  A (r) x + F,
\end{equation}
where
\[
A(r)= - \mat{h_1 \\ h_3} [h_{1,r}, h_{3,r}] = \frac{m}{r}h_1\mat{
h_{1}h_3,& -h^2_1\\h^2_3,&-h_{1}h_3}
\]
and
\[
F = \mat{F_1 \\ F_3}, \quad
F_j = - (\de \bxi_r \cdot \be_a)v_{a,j} - (\bxi_{b,r} \cdot \de \be)h_j
 - (\bv_{b,r} \cdot \de \be)\xi_{a,j}- (\bv_{b,r}
\cdot \be_b)\de \xi_j, \quad j=1,3.
\]

To simplify the linear part $\tilde x_r = A(r)\tilde x$, let
$y=U^{-1}\tilde x$ where
\[
U(r)=\begin{bmatrix} h_{1},& -h_3\\h_3,& h_{1}\end{bmatrix},\quad
U^{-1}=\begin{bmatrix} h_{1},& h_3\\-h_3,& h_{1}\end{bmatrix}.
\]
Then $y$ satisfies
\[
y_r= (U^{-1})_r  U y + U^{-1} A U y
=\frac{m}{r}h_1\begin{bmatrix} 0,& 0\\-1,& 0\end{bmatrix}y.
\]
This linear system can be solved explicitly,
\[
y(r) =\mat{ 1,& 0\\ p(\rho,r), & 1} y(\rho),
\quad p(\rho,r)= -\bke{\int_{\rho}^r\frac{m}{r}h_1(\tau)d\tau} =
-\big [ 2 \arctan \tau^m \big ]_\rho^r.
\]
Thus the linear system $\tilde x_r = A(r)\tilde x$ has the
solution $\tilde x(r)=P(\rho,r) \tilde x(\rho)$ with the propagator
\[
P(\rho,r) = U(r)\mat{ 1,& 0\\ p(\rho,r), & 1}U^{-1}(\rho).
\]
The original system \eqref{lwp140} with $x(0)=0$ has the solution
\begin{equation*}
x(r)=\int^r_0 P(\rho,r)F(\rho)d\rho.
\end{equation*}

To estimate $x(r)$, the two terms of $F_3$ with $h_3$ as the last factor,
\[
\tilde F_3 = - (\de \bxi_r \cdot \be_a)h_3 - ( \bxi_{b,r} \cdot \de \be) h_3
\]
require special care since it may not be in $L^1(dr)$.
Other terms can be estimated as follows:
\begin{equation*}
\Big |\int_0^r P(\rho,r)\mat{F_1\\ F_3 -\tilde F_3}d\rho \Big|
\lec \int_0^\infty |F_1| + |F_3- \tilde F_3| dr
\lec
\norm{\delta z}_{X} +\bke{\norm{z_a}_X+\|z_b \|_{X}} \norm{\delta
\be}_{L^{\infty}}.
\end{equation*}
We treat $\tilde F_3 $ by integration by parts:
\begin{multline*}
\int^r_0  P(\rho,r)  \mat{ 0\\ \tilde F_3} d\rho
= \int^r_0 P(\rho,r) \mat{ 0\\ -(\delta\bxi_{\rho}\cdot \be_a
+\bxi_{b,\rho}\cdot \de \be)h_3} d\rho
\\
 = - \mat{0 \\ (\delta\bxi\cdot \be_a +\bxi_b \cdot \de \be)h_3}(r)
+\int^r_0 P(\rho,r) \mat{ 0\\ (\de \bxi\cdot \be_{a,\rho}
+ \bxi_b\cdot \de e_\rho)h_3 + (\delta\bxi \cdot \be_a
+\bxi_b \cdot \de \be)h_{3,\rho}} d\rho
\\
  +\int^r_0 P_{\rho}(\rho,r)\mat{0 \\
(\delta\bxi\cdot \be_a)h_3 +(\bxi_b\cdot\de \be)h_3}d\rho = \sum _{j=1}^3
I_j .
\end{multline*}
Now we estimate the right side one by one. For $I_1$,
\[
|I_1| \lec
\norm{\de \bxi}_{L^{\infty}}+
\|\bxi_b \|_{L^{\infty}}
\|\de \be \|_{L^{\infty}}
\lec \norm{\de z}_X+\|z_b \|_X\norm{\de \be}_{L^{\infty}}.
\]
For $I_2$, observe that
\[
\norm{\be_{a,r}}_{L^2}\leq C, \quad \norm{\de \be_r}_{L^2}
\lec \norm{\de z}_X+\norm{\de \be}_{L^{\infty}},
\]
due to the facts that $\be_{a,r}=-(\bv_{a,r}\cdot \be_a)\bv_a$ and
$\de \be_r=-(\bv_{a,r}\cdot \be_a)\bv_a + (\bv_{b,r}\cdot
\be_b)\bv_b$.
Thus
\[
|I_2| \lec \norm{\de z}_X+ \norm{z_b}_X \bke{\norm{\de z}_X+\norm{\de
\be}_{L^{\infty}}}.
\]
To estimate the last term $I_3$, note  that
\[
P_{\rho}(\rho,r)
=\frac{m}{\rho}h_1(\rho)
U(r) \cdot
\bket{\begin{bmatrix} 0& 0\\1& 0\end{bmatrix}\cdot
U^{-1}(\rho)
+\begin{bmatrix} 1& 0\\p(\rho,r)
& 1\end{bmatrix}\cdot
\begin{bmatrix} -h_{3},& h_1\\-h_1,& -h_{3}\end{bmatrix}(\rho)}
\]
and hence $|P_{\rho}(\rho,r)| \lec h_1(\rho)/\rho$. We get
\[
|I_3| \lec
\int^r_0 \abs{\frac{h_1}{\rho}}\bke{\abs{(\delta\bxi\cdot \be_a)h_3}
+\abs{(\bxi_b\cdot\de \be)h_3}}d\rho
\lec \|\frac{h_1}\rho \|_{L^2}
(\norm{\de z}_X+\|z_b\|_X\norm{\de \be}_{L^{\infty}}).
\]

Summing up, we have shown
\[
\norm{\de \be}_{L^{\infty}}\lec \norm{\de z}_X
+\bke{\norm{z_a}_X+\| z_b \|_X}\norm{\de \be}_{L^{\infty}}.
\]
Since $\norm{z_a}_X+\norm{z_b}_X \ll 1$, we can absorb the last term
to the left side. The lemma is proved.
\end{proof}

\begin{lemma} \label{lwp:T3}
For given $s>0$, $\al \in \R$, and $q\in L^2_{rad}$ with
$\norm{q}_{L^2} \le \de$, there is a unique function $z=Z(q,s,\al) \in
X$ so that $\bka{h_1,z}_X=0$, $\|z\|_X \lec \de$ and the function
$\bv=\BV(Z,s,\al)$ satisfies \eqref{lwp-10}. Moreover, $Z(q,s,\al)$ is
independent of $\al$ and continuous in $q$ and $s$.
\end{lemma}

\begin{proof}
Simple comparison shows
\begin{equation}
Z(q,s,\al) = Z(q(\cdot s), 1, 0).
\label{lwp-13}
\end{equation}
Thus it suffices to prove the case $s=1$ and $\alpha=0$. We will
construct $Z(q,1,0)$ by a contraction mapping argument. Define the map
\begin{equation}
\Phi^q (z)(r) =  L_0^{-1} \Pi [q \BE(z) - {\bf G_0}(z)](r),
\end{equation}
where $\Pi = (\bV^r)^{-1} P^{h(r)}$ is a projection of vector fields on
$\R^+$ to $L^2(rdr)$, with the mapping $(\bV^r)^{-1}:T_{\bh(r)}\S^2 \to \C
$ and the projection $P^{\bh(r)}:\R^3 \to T_{\bh(r)}\S^2$ defined in
Section \ref{sec:z}; $ L_0^{-1}$ is the inverse map of $L_0$ and maps
$L^2(rdr)$ to the $X$-subspace $h_1^\perp$; $\BE(z)$ is defined after
\eqref{lwp110}, and ${\bf G_0}(z)$ is defined by \eqref{lwp-15}.

We will show that $\Phi^q$ is a contraction mapping in the class
\[
{\cal A}_\de= \bket{ z \in X: \norm{z}_X
\le 2   C_1 \de}, \quad C_1 = \| L_0^{-1} (\bV^r)^{-1}
P^{\bh(r)} \|_{B(L^2,X)}
\]
for sufficiently small $\de>0$. First,
\[
\norm{\Phi^q(z)}_X \le C_1 \norm{q}_2 + C \norm{{\bf G_0}}_2 \le C_1 \de
 + C\norm{z}_X^2.
\]
Thus $\Phi^q$ maps $\cal A$ into itself if $\de$ is sufficiently
small.  We now prove difference estimates for $\Phi^q$. Suppose
$z_a,z_b \in {\cal A}$ are given and let $\bv_l=\BV(z_l)$ and $\be_l =
\BE(z_l)$, $l=a,b$.  Also define $\bxi_l$ by \eqref{lwp-05} and note
$\de \bxi = \de \bv$.  By Lemma \ref{lwp:T2},
\[
\norm{\de \bv}_{X}+\norm{\de \bxi}_{X}+\norm{\de \be}_{L^\infty} \lec
\norm{\de z}_X.
\]
We now estimate $\de {\bf G_0}(z) = {\bf G_0}(z_a)-{\bf G_0}(z_b)$ in
$L_p$, $p=2,4$ (we need $p=4$ later):
\begin{equation}\label{lwp-45}
\norm{\de {\bf G_0}(z)}_{L_p}\lec \norm{\delta\gamma_r}_{L_p}
+\nrm{\frac{\delta\gamma}{r}}_{L_p} + \norm{\de(\xi_3 \bxi)}_{L_p}
\lec (\norm{z_a}_X + \|z_b\|_X) \norm{\de z}_{X_p}.
\end{equation}
Thus,
\begin{equation}\label{lwp-50}
\norm{\Phi^q(z_a) - \Phi^q(z_b)}_X \lec \norm{q \de \be - \de {\bf G_0}(z)}_{L^2}
\lec \norm{q}_{L^2} \norm{\de \be}_{L^\infty} + (\norm{z_a}_X +
\|z_b\|_X) \norm{\de z}_X \ll \norm{\de z}_X.
\end{equation}
Thus $\Phi^q$ is indeed a contraction mapping and the function
$Z(q,s,\al)$ exists.

We now consider the continuity. The continuity in $s$ follows from
\eqref{lwp-13}, although it may not be H\"older continuous.  For the
continuity in $q$, let $q_a$ and $q_b$ be given and $z_l =
Z(q_l,s,\al)$, $l=a,b$. An estimate similar to \eqref{lwp-50} shows
\begin{equation}
\norm{\de z}_X = \norm{\Phi^{q_a}(z_a) - \Phi^{q_b}(z_b)}_X \lec
\norm{\de q}_{L^2} + \e \norm{\de z}_X ,
\end{equation}
where $\e= \norm{q_a}_{L^2} + \| q_b\|_{L^2}+\norm{z_a}_X + \|z_b\|_X
\ll 1$ and hence $\e \norm{\de z}_X$ can be absorbed to the left side.
This shows continuity in $q$ in $L^2$-norm.
%
\end{proof}


\subsection{Evolution system of $q$, $s$ and $\al$}
\label{lwp:S2}

By \eqref{lwp-05}, the dynamics of $\bu$ is completely determined by the
dynamics of $z$, $s$ and $\al$.  Because of Lemma \ref{lwp:T3}, it is
also completely determined by the dynamics of $q$, $s$ and $\al$. The
latter system is preferred by us since the $q$ equation is easier than
the $z$ equation to estimate, and $q$ lies in a more familiar space
$L^2$, rather than $z$ in $X$.

The equations for $z$ and $q$ are given by \eqref{eq:z} and
\eqref{eq:q}, respectively.  However, since we choose the
orthogonality condition \eqref{lwp:orth}, i.e., $\bka{h_1,z(t)}=0$ for
all $t$, the equations for $s$ and $\al$ are different from
\eqref{eq:ode}.

We now specify the equations we will use.  Let $\tilde{q}
:=e^{i(m+1)\theta}q$.  Recall $\nu \be = \nu_1 \be+ \nu_2 J^v \be =
J^{\bv}R \bv = \hk - v_3 \bv $ and $\nu_r = v_3(q+\frac mr \nu)$.  By
\eqref{eq:q} and an integration by parts on the potential defined in
\eqref{eq:NL}, we obtain
\begin{equation}
i\tilde{q}_t+\Delta\tilde{q}=V \tilde q, \quad
V = V_1 - V_2 + \int_r^\infty \frac 2{r'}V_2(r') dr'
\label{lwq210}
\end{equation}
where
\[
V_1:= \frac{m(1+v_3)(mv_3-m-2)}{r^2}+ \frac{mv_{3,r}}{r},
\quad
V_2 := \frac 12 |q|^2 + \Re \frac mr \bar \nu q.
\]

For $s$ and $\al$, condition \eqref{lwp:orth} implies $\bka{h_1, \pd_t
z(t)}_X = 0$. Substituting in \eqref{eq:z}, we get
\begin{equation}\label{lwp215}
 \bka{ h_1,\, \bke{s^2\dot \al - im s \dot s } (1+\gamma)h_1 +
 s^2\dot \al i z h_3 - s \dot s r z_r}_X
= \bka{h_1, \ -i Nz + P{\bf F_1} }_X.
\end{equation}
Note
\[
\bka{h_1, \  Nz}_X = ( L_0 N_0h_1,L_0z)_{L^2}, \quad
\bka{h_1, \ r\pd_r z}_X = ( rN_0 h_1,z_r)_{L^2}.
\]
Let $G_1:= \bka{h_1,P{\bf F_1}}_X = (N_0 h_1,P{\bf F_1})_{L^2}$
where $N_0 := - \Delta _r + \frac {m^2}{r^2}$.
By Lemma
\ref{Ap2:T2} with $g=N_0 h_1$,
\begin{align*} 
G_1 =
 \int_0^\infty \bigg( & i g_r (- \gamma z_r + z \gamma_r)
+ \frac mrh_1 g (-2 \gamma_r - i z_2 z_{1,r} + i z_1 z_{2,r})
\\
& + \frac mr(h_1g)_r (\gamma^2 - i z_2 z)
+ i \frac {m^2}{r^2}(2h_1^2-1) g \gamma z \bigg) rdr.
\end{align*}
Separating real and imaginary parts, we can rewrite \eqref{lwp215}
as a system for $\dot \al $ and $\dot s$:
\begin{equation}\label{lwp220}
\bke{\norm{h_1}^2_X I + A} \mat{s^2 \dot \al \\ -m s \dot s}
= \vec {G}_2 :=
%
 \mat{ (L_0N_0h_1,L_0z_2)_2 \\ -(L_0N_0h_1,L_0z_1)_2}
+ \mat{ \Re G_1 \\ \Im G_1},
\end{equation}
where
\[
I=\begin{bmatrix} 1,& 0
\\0,& 1
\end{bmatrix},\quad
A=\begin{bmatrix} \bka{ h_1,\, \gamma h_1- z_2 h_3}_X,
& \frac 1m (rN_0 h_1,\, z_{1,r})_{L^2}
\\
\bka{h_1,\, z_1 h_3}_X,& \bka{ h_1,\, \gamma h_1 }_X
+\frac 1m (rN_0 h_1, z_{2,r})_{L^2}
\end{bmatrix}.
\]
We have $\norm{A}_{L^\infty} \lec \norm{z}_X$.

We will study the integral equation version of \eqref{lwq210} and
\eqref{lwp220} for $\tilde q$, $s$, and $\al$:
\begin{equation}
\tilde q(t) = e^{-it \Delta} \tilde q_0 -i \int_0^t e^{-i(t-\tau) \Delta}
(V\tilde q)(\tau) \, d\tau,
\label{lwp230}
\end{equation}
\begin{equation}
\mat{ s(t)\\ \al(t)} = \mat{ s_0 \\ \al_0} + \int _0^t
\bket{ \mat{0 &-
(ms)^{-1} \\ s^{-2} & 0} \bke{\norm{h_1}^2_X I + A} ^{-1}
\vec {G}_2}(\tau) \,d\tau.
\label{lwp240}
\end{equation}


\subsection{Contraction mapping and conclusion}
\label{lwp:S3}

Let $q_0 \in L^2_{rad}(\R^2)$, $s_0>0$, and $\al_0 \in \R$ be given,
with $\norm{q_0}_{L^2} \le \de$. For $\de,\si>0$ sufficiently small we
will find a solution of \eqref{lwp230}--\eqref{lwp240} for $ t\in
I=[0,\si s_0^2]$.

We will first construct the solution assuming
$s_0=1$.  The solution for general $s_0$ is obtained from rescaling,
\[
  \bu(t,x) = \tilde \bu(t/s_0^2, x/s_0)
\]
where $\tilde u$ is the solution corresponding to initial data $\tilde
\bu_0(x) = \bu_0(x/s_0)$, and $s(\tilde \bu_0) = 1$.

Assuming $s_0=1$, we will define a (contraction) mapping in the
following class
\begin{multline}
{\cal A}_{\de,\si} = \big \{ (\tilde q,s,\al):
I=[0,\si] \to L^2(\R^2) \times \R^+ \times \R:
\\
\norm{\tilde q}_{Str[I]}\le  \de; \   \forall t, \
q(t)= e^{-i(m+1)\th}\tilde q(t) \in L^2_{rad},
\ s(t) \in [0.5,1.5]
\big \},
\end{multline}
for sufficiently small $\de,\si>0$.
Here
\[
\norm{q}_{Str[I]} \equiv
\norm{q}_{L^\infty_t L^2_x [I]\cap L^4_t L^4_x[I]
\cap L^{8/3}_t L^8_x [I]
}.
\]
The map is defined as follows. Let $\tilde q_0 = e^{i(m+1)\th}q_0$.
Suppose $Q= (\tilde q,s,\al)(t) \in {\cal A}_{\de,\si}$ has been
chosen. For each $t \in I$, let $ q = e^{-i(m+1)\th} \tilde q$, let
$z=Z(q,s,\al)$ be defined by Lemma \ref{lwp:T3}, and let
$\bv=\BV(z,s,\al)$ and $\be=\BE(z, s, \al)$ be defined by \eqref{lwp-05}
and \eqref{lwp110}, respectively. We then substitute these functions
into the right sides of \eqref{lwp230} and \eqref{lwp240}. The output
functions are denoted as $\tilde q^\sharp(t)$, $s^\sharp(t)$, and
$\al^\sharp(t)$. The map $Q \to \Psi(Q)= ( \tilde q^\sharp, s^\sharp,
\al^\sharp$) is the (contraction) mapping.

The following estimates are shown in \cite[Lem.~3.1]{GKT}.
\begin{equation}
\label{lwp320}
\|\tilde q^\sharp \|_{Str[I]} \lec
\norm{q_0}_{L^{2}_{x}}+
(\sigma^{\frac{1}{2}} +\norm{q}^2_{L^4_{t,x}[I]})
\norm{q}_{L^4_{t,x}[I]}.
\end{equation}
We also have $ |\vec{G}_2| \lec \|z\|_X + \|z\|_X^4$ and thus
\[
|s^\sharp (t)-1| + |\al^\sharp(t)-\al_0|
\lec \int_0^t  |\vec{G}_2(\tau)| \,d\tau
\lec \si \norm{q}_{L^\infty_t L^2_x} + \si\norm{q}_{L^\infty_t L^2_x}^4.
\]
Therefore ${\cal A}_{\de,\si}$ is invariant under the map $\Psi$ if
$\de$ and $\si$ are sufficiently small.

\medskip

We now consider the more delicate difference estimate. Suppose we
have $Q_l=(\tilde q_l,s_l,\al_l)(t)$ for $l=a,b$. Let $z_l$, $v_l$,
$\be_l$, $\tilde q_l^\sharp$, $s_l^\sharp$ and $\al_l^\sharp$ be
defined respectively. Denote
\begin{equation} \label{lwp335}
\de \tilde q = \tilde q_a(t,r) - \tilde q_b(t,r), \quad
\de z = z_a(t,r/s_a)-z_b(t,r/s_b), \quad etc.
\end{equation}
Note that we define $\de z$ in terms of $r$, not in $\rho$, i.e., $\de
z \not = z_a(\rho)-z_b(\rho)$.  See Remark \ref{lwp:rk1} after the proof.  In
the rest of the proof, we denote
\[
\norm{q}_2 = \max_{a,b} (\norm{q_a}_2, \norm{q_b}_2 ),\quad
\norm{z}_X = \max_{a,b} (\norm{z_a}_X, \norm{z_b}_X), \quad etc.
\]

To start with, note that
\begin{equation}\label{lwp336}
\norm{z}_{L^\infty_t X} \lec \de, \quad
|\de h_1| \lec |\de s| \frac {h_1}r, \quad
|\de h_3| \lec |\de s| \frac {h_1^2}r, \quad
|\de \gamma| \lec |z||\de z|.
\end{equation}

We first estimate $\de \be= \be_a-\be_b= \BE(z_a,s_a,\al_a) -
\BE(z_b,s_b,\al_b)$.  By \eqref{lwp115},
\[
|\de  \be | \lec |\de \al|
 + \nrm{\BE (z_a(\cdot/s_a)) - \BE (z_b(\cdot/s_b)}_{L^\infty}.
\]
By Lemma \ref{lwp:T2}, $\nrm{\BE (z_a(\cdot/s_a)) - \BE
(z_b(\cdot/s_b)}_{L^\infty} \lec
\nrm{z_a(\cdot/s_a) - z_b(\cdot/s_b)}_X =
\nrm{\de z}_X$. Thus
\begin{equation}\label{lwp338}
|\de  \be | \lec |\de \al| +\nrm{\de z}_X.
\end{equation}

We next estimate $\norm{\de z}_{X}$. By \eqref{lwp-10},
\[
\de (L_0 z)\hj = \de [s e^{-\al R}q \be(r) ] - \de {\bf G_0} .
\]
Here $\de (L_0 z) = L_0(\frac r{s_a}) z_a(\frac r{s_a})
- L_0(\frac r{s_b}) z_b(\frac r{s_b})$ and
$\de {\bf G_0} = {\bf G_0}(z_a(\frac r{s_a})) - {\bf G_0}(z_b(\frac r{s_b}))$.
Rewrite
\[
\de (L_0 z) = D_1 + L_0(r/s_a) \de_1 z,
\]
where
\[
D_1 = (\de L_0) z_b(r/s_b) ,\quad
\de_1 z= z_a(r/s_a) - \Pi_{s_a} z_b(r/s_b)
\]
and $\Pi_s$ is the projection removing $h_1(z/s)$: $\Pi_s f = f -
\frac{\bka{h_1(\cdot/s), f}_X} {\bka{h_1, h_1}_X}h_1(\cdot/s)$. Here
we have used $L_0(r/s_a) = L_0(r/s_a) \Pi_{s_a}$. Since $L_0(r/s) =
s[\pd_r - \frac mr h_3(r/s)]$, we have $ \de L_0 \sim \de s [L_0(r/s)
- s \frac {m^2}{r^2} h_1^2(r/s)\cdot \frac r{s^2}]$, and hence
\[
\norm{D_1}_{L^2} \lec |\de s| \cdot \norm{z}_{X}.
\]

Thus, taking $L_0(r/s_a)^{-1}$,
\[
\norm{\de_1 z}_{X}
\lec \norm{\de [s e^{-\al
R}q \be(r) ]}_2 + \norm{\de {\bf G_0}}_2 + \norm{D_1}_{L^2}.
\]

We can decompose
\[
\de z = \de_1 z + \de_2 z, \quad
\de_2 z = (1-\Pi_{s_a}) z_b(r/s_b),
\]
and we have
\[
\norm{\de_2 z}_{X} \lec
\bka{h_1(z/s_a) - h_1(z/s_b),  z_b(r/s_b)}_{X}
\le |\de s| \norm{z}_{X}.
\]

Note
\begin{align*}
|\de {\bf G_0}| &\lec |\de {\bf h}| |\gamma_\rho| + \frac {\de s}r(|\gamma|+|\bxi|^2)
+ |\de \gamma_\rho| + \frac 1r(|\de \ga| + |\bxi||\de \bxi|)
\\
&\lec |\de s|(|z||z_r| + |z|^2/r) + |\de z|(|z_r|+|z|/r) + |z| |\de z_r|.
\end{align*}
Thus,
\[
\norm{\de {\bf G_0}}_2 \lec |\de s| \norm{z}_{X}^2 + \norm{z}_{X} \norm{\de
z}_{X} .
\]

Finally,
\[
\norm{\de [s e^{-\al R}q \be(r) ]}_2 \lec (|\de s| + |\de \al| +
\norm{\de \be}_{L^\infty}) \cdot \norm{q}_2 + \norm{\de q}_2.
\]

Adding these estimates, using \eqref{lwp338}, and $\norm{z}_{X}
\lec \norm{q}_2$,  we get
\[
\norm{\de z}_{X} \lec (|\de s| + |\de \al| + \norm{\de z}_{X})
\cdot \norm{q}_2 + \norm{\de q}_2.
\]
Absorbing $\norm{\de z}_{X} \norm{q}_2$ to the left side, we get
\begin{equation} \label{lwp345}
\norm{\de z}_{X} \lec (|\de s| + |\de \al|) \cdot
\norm{q}_2 + \norm{\de q}_2.
\end{equation}

We now estimate $\norm{\de \tilde q^\sharp}_{Str[I]}$. Apply Strichartz
estimate to the difference of \eqref{lwp230},
\begin{align*}
\norm{\de \tilde q^\sharp}_{Str[I]}
&\lec \norm{\de (V\tilde
q)}_{L^{4/3}_{t,x}} \lec
\norm{V(\de \tilde q)}_{L^{4/3}_{t,x}}
+ \norm{(\de V)\tilde q}_{L^{4/3}_{t,x}}.
\\
&\lec \norm{V}_{L^{2}_{t,x}}\norm{\de \tilde q}_{L^{4}_{t,x}}
+ \norm{\de V}_{L^2_{t,x} + L^{8/3}_tL^{8/5}_x}
\norm{\tilde q}_{L^{4}_{t,x}\cap L^{8/3}_tL ^8_x}.
\end{align*}
Recall  $V=V_1-V_2 + \int^r \frac 1{r'} V_2$. By the 4-dimensional
Hardy inequality, for each fixed $t$,
\[
\norm{V}_{L^2_x} \lec \norm{V_1}_{L^2_x} + \norm{V_2}_{L^2_x}
\lec \nrm{\frac {1+v_3}{r^2}}_2 +
 \nrm{\frac {v_{3,r}}{r}}_2
+ \norm{q}_4^2 + \norm{q}_4 \cdot \nrm{\frac \nu r}_4,
\]
and, since $v_3(r) = (h_3+ h_3 \ga + h_1 z_2)(r/s)$ and $|\nu|= |\hk -
v_3 \bv|$,
\[
\nrm{\frac {1+v_3}{r^2}}_2 + \nrm{\frac {v_{3,r}}{r}}_2
\lec 1+ \norm{z}_X+ \nrm{\frac zr}_4\cdot \nrm{z}_{X_4}  ,
\]
\[
\nrm{\frac \nu r}_4^2 = \nrm{\frac {1-v_3^2}{ r^2}}_2
\lec \nrm{\frac {1+v_3}{r^2}}_2.
\]
Thus $\norm{V}_{L^2_x} \lec 1 + \norm{q}_{L^4_x}^2$ and
hence $\norm{V}_{L^2_{t,x}[I]} \lec \si^{1/2}+
\norm{q}_{L^4_{t,x}}^2$.

Denote $Y=L^2_{t,x} + L^{8/3}_tL^{8/5}_x$.  By Hardy inequality again,
\begin{align*}
\norm{\de V}_Y &\lec \norm{\de V_1}_Y + \norm{\de V_2}_Y \\
&\lec
\nrm{\frac {\de v_3}{r^2}}_Y + \nrm{\frac {\pd_r \de v_{3}}{r}}_Y +
(\norm{q}_{L^4_{t,x}}+\nrm{\frac \nu r}_{L^4_{t,x}})\norm{\de
q}_{L^4_{t,x}} + \norm{q}_{L^{8/3}_tL^{8}_x} \cdot \nrm{\frac {\de
\nu} r}_{L^{\infty}_tL^{2}_x}.
\end{align*}

Note $\nu = \be \cdot (\hk - v_3 \bv)$.
Thus $\de \nu = \de \be \cdot
(\hk - v_3 \bv) - \be \cdot ((\de v_3) \bv + v_3 \de \bv)$, and
\[
\nrm{\frac {\de \nu}r}_{L^2_x} \lec \nrm{\de \be}_\infty \nrm{\frac 1r(\hk -
v_3 \bv)}_2 + \nrm{\frac 1r \de \bv}_2.
\]
Since $\nrm{\frac 1r(\hk - v_3 \bv)}_2 \lec 1 + \norm{z}_X^2 \lec 1 $ and
$\nrm{\frac 1r \de \bv}_2 \lec |\de \al|\nrm{\frac{\bh+\bxi}r}_2 +
\nrm{\frac{\de \bh}r}_2+ \nrm{\frac{\de z}r}_2 $,
we conclude using \eqref{lwp338} and \eqref{lwp345},
\[
\nrm{\frac {\de \nu}r}_{L^2_x} \lec |\de s| + |\de \al|+ \nrm{\de
  q}_2.
\]

For $\frac {\de v_3}{r^2}$ and $\frac {\pd_r \de v_{3}}{r}$,
since $v_3(r) = (h_3+ h_3 \ga + h_1 z_2)(r/s)$,
\[
\frac 1{r^2}\, |\de v_3| \lec \frac 1{r^2}\, (|\de h_3| + |\de h_1|
|z| + |\de \ga| + h_1 |\de z|) \lec
\frac{h_1+|z|}r \bke{ |\de s|\,\frac {h_1}r +\frac{|\de z|}r},
\]
\[
\frac 1r \, |\pd_r \de v_{3}| \lec
|\de s|\, \bke{ \frac {h_1(h_1+|z|)}{r^2} + \frac {h_1 + h_1^2|z|}r |z_r|}
+\frac {h_1 + h_1^2|z|}r  \frac {|\de z|}r
+ \frac {h_1 + |z|}r  |\pd_r \de z|.
\]
We do not want to bound $\frac zr \frac {\de z}r$ and $\frac zr \pd_r
\de z$ in $L^2_x$ since otherwise we would need a bound for $\norm{\de
z}_{X_p}$, $p>2$, which requires extra effort.  We have
\begin{align*}
\nrm{\frac {\de v_3}{r^2}}_Y + \nrm{\frac {\pd_r \de v_{3}}{r}}_Y
&\lec \norm{\de s}_{L^\infty_t}
\norm{\bke{ \frac {h_1(h_1+|z|)}{r^2} + \frac {h_1 +
h_1^2|z|}r |z_r|} }_{L^2_{t,x}}
\\
& \quad
+\norm{\frac {h_1 + h_1^2|z|}r \frac {|\de z|}r }_{L^2_{t,x}}
+ \norm{\frac {h_1 + |z|}r ( \frac {| \de z|}r+|\pd_r \de z|)}
_{L^{8/3}_tL^{8/5}_x}
\\
&\lec \norm{\de s}_{L^\infty_t} +
(1+ \norm{z/r}_{L^{8/3}_tL^{8}_x})\norm{\de z}_{L^\infty_t X} .
\end{align*}

Using $\norm{z/r}_{L^{8/3}_tL^{8}_x} \lec
\norm{q}_{L^{8/3}_tL^{8}_x}\lec \de$, and
\eqref{lwp345}, we conclude
\begin{equation} \label{lwp360}
\nrm{\de \tilde q^\sharp}_{Str[I]}
\lec
(\si^{1/2}+\norm{\tilde q}_{L^{4}_{t,x}}^2)\norm{\de \tilde q}_{L^{4}_{t,x}}
+ \norm{\tilde q}_{Str[I]}
 (\norm{\de s}_{L^\infty_t} +\norm{\de \al}_{L^\infty_t}
+ \norm{\de \tilde q}_{L^{\infty}_{t}L^2_x} ).
\end{equation}

We now estimate $\de s^\sharp$ and $\de \al^\sharp$. Estimating the
difference of \eqref{lwp240},
\[
\nrm{\de s^\sharp}_{L^\infty(I)} +\nrm{\de \al^\sharp}_{L^\infty(I)}
\lec \int_I (|\de s| + |\de A|) |\vec{G}_2| + |\de \vec{G}_2| d \tau.
\]
Note that $|\vec{G}_2| \lec \norm{z}_X + \norm{z}_X^4$,
\[
|\de A| \lec \nrm{\de \bh}_X \nrm{z}_X + \nrm{ h_1}_X \nrm{\de z}_X
\lec |\de s| \nrm{z}_X + \nrm{\de z}_X ,
\]
and
\begin{align*}
|\de \vec{G}_2| &\lec \nrm{\de h}_X \nrm{z}_X + \nrm{ h_1}_X \nrm{\de z}_X
+|\de G_1|
\\
& \lec|\de s| \nrm{z}_X + \nrm{\de z}_X +
(1+\nrm{z}_\infty)\bke{\nrm{z}_\infty \nrm{\pd_r \de z}_2 + \nrm{\pd_r  z}_2
\nrm{\de z}_\infty}
\\
&\quad + (\nrm{z}_\infty+\nrm{z}_\infty^3) \nrm{\de z/r}_2,
\end{align*}
Thus,
\begin{align}
\nonumber
\nrm{\de s^\sharp}_{L^\infty(I)} +\nrm{\de \al^\sharp}_{L^\infty(I)}
&\lec
\int_I |\de s| \norm{ z}_X + (1+\norm{ z}_X^3)\norm{\de z}_X d\tau
\\
&\lec
\si \norm{ z}_X \nrm{\de s}_{L^\infty(I)}
+ \si \norm{\de z}_{L^\infty_t X} .
\label{lwp380}
\end{align}

Combining \eqref{lwp345}, \eqref{lwp360} and \eqref{lwp380}, we have
proved that
\begin{equation} \label{lwp390}
\nrm{\de \tilde q^\sharp }_{Str[I]}
+ \nrm{\de s^\sharp}_{L^\infty(I)} +\nrm{\de \al^\sharp}_{L^\infty(I)}
\lec (\sigma^{1/2} + \de)
\bke{\nrm{\de \tilde q}_{Str[I]}
+ \nrm{\de s}_{L^\infty(I)} +\nrm{\de \al}_{L^\infty(I)} }.
\end{equation}
Thus $\Psi$ is a contraction mapping on ${\cal A}_{\de,\si}$ if $\si$
and $\de$ are sufficiently small. 
We have therefore established the unique existence of a triplet 
$[s_W(t), \alpha_W(t), q_W(t)]$ solving the $(s,\alpha,q)$-system.
This yields a map $\bu_W(t) \in C([0,T];\Sigma_m)$.

If $\bu_0 \in \dot H^2$, the a priori estimates in 
\cite[Lem.~3.1]{GKT} show
$\norm{\nabla \tilde q}_{Str[I]}$ is uniformly bounded,
so $\bu_W(t) \in C(I;\Sigma_m \cap \dot H^2)$. 

If $\bu_0^n \to \bu_0$ in $\Sigma_m \cap \dot H^k$, $k=1,2$,
a difference estimate similar to \eqref{lwp390} shows
\[
  D^n \lec \norm{\tilde q^n_0 - \tilde q_0}_2
  +  (\sigma^{1/2} + \de)D^n.
\]
where $D^n = \nrm{\tilde q^n- \tilde q}_{Str[I]} + \nrm{s^n-
s}_{L^\infty(I)} +\nrm{\al^n - \al}_{L^\infty(I)} $.  Thus $D_n \to 0$
as $n \to \infty$, and hence $\bu^n_W \to \bu_W$.

The energy $\sE(\bu_W(t))$ is conserved since 
$\sE(\bu_W(t)) = 4 \pi m + \pi
\norm{q(t)}_{L^2_x}^2 = 4 \pi m + \pi \norm{q_0}_{L^2_x}^2$. 

Finally, we must verify that $\bu_W$ is a solution of the 
Schr\"odinger flow as in Definition~\ref{def:sol}. 
To do this, approximate the initial data $\bu_0$
in $\Sigma_m$ by $\bu_0^k$ with $\nabla \bu_0^k \in H^{10}$ (say).
By~\cite{SSB} there is a unique strong solution
$\bu_S^k(t)$ with initial data $\bu_0^k$.
The corresponding triple $[s^k_S(t), \alpha^k_S(t), q^k_S(t)]$ 
must satisfy the $(s,\alpha,q)$-system.
By uniqueness, $s^k_S(t) \equiv s^k_W(t)$, etc.,
and so $\bu_W^k(t) \equiv \bu_S^k(t)$.
By continuous dependence on $\dot H^1$ data, 
$\bu_S^k$ converges to $\bu_W$ in
$C([0,T];\Sigma_m)$, and in particular in $C([0,T];L^2_{loc})$.
Finally, $\bu_S^k$ satisfies the weak form of the
Schr\"odinger flow (Definition~\ref{def:sol}), and passing
to the limit, so does $\bu_W$.
Dropping the subscript $W$ ($\bu := \bu_W$), 
Theorem \ref{lwp:T1} is established.

We now consider Corollary \ref{lwp:T1b}. 
Suppose $T$ is the blow-up time. By Theorem~\ref{lwp:T1}, 
for each $t <T $ we have $T-t \ge \si s(\bu(t))^2$.  
Thus $s(\bu(t)) \le \si^{-1/2} \sqrt{T-t}$. 
If $k=2$, by \cite[Th.~2.1]{GKT}, 
$\norm{\bu(t)}_{\dot H^2} \ge C_2 / s(\bu(t)) \ge C_2
\si^{1/2} (T-t)^{-1/2}$. On the other hand, the 
$\dot H^2$-estimates of \cite{GKT} show that
the $\dot H^2$-norm can only blow-up if 
$\liminf_{t \to T^-} s(t) = 0$. Thus $T^2_{max} = T^1_{max}$.
Statement (ii) follows from Theorem
\ref{lwp:T1} directly.  Corollary \ref{lwp:T1b} is established.

\begin{remark}\label{lwp:rk1}
\begin{enumerate}

\item In Theorem \ref{lwp:T1}, we did not try to prove continuity
on data $\bu_0$ in $\dot H^2$, which would require difference
estimates in $H^1$ for $\tilde q$.

\item In \eqref{lwp335}, we define $\de z$ in terms of $r$, not in
$\rho$, i.e., $\de z \not = \tilde \de z= z_a(\rho)-z_b(\rho)$.
Indeed, in view of \eqref{lwp-10}, since $L_0$ depends on $\rho$, it
may seem natural to bound $\tilde \de z$ using $L_0 \tilde \de z \hj =
\de [s e^{-\al R}q \be(s \rho) ] + \tilde \de {\bf G_0}$.  However, to bound
the right side we need to bound the difference $ q_a\be_a(s_b\rho)
-q_a\be_a (s_a\rho) = \int_{s_a}^{s_b} \rho \pd_r(q_a\be_a)(\sigma
\rho) d \sigma $, for which $\norm{\bu}_{\dot H^2}$ is insufficient and
we need a weighted norm of $\bu$.  The reason is that the dilation
magnifies the difference when $\rho$ is large.  In addition, to bound
$\de v_3$ using $\tilde \de z$ instead of $\de z$, one needs a bound
for $z_{rr}$.

\item In the proof we have avoided using $\norm{\de z}_{X_4}$
since its estimate requires $\norm{\de \be}_{\infty}$. We know how
to control $\norm{\de \be}_\infty$ by  $\norm{\de z}_{X}$, but we
do not know if $\norm{\de \be}_{\infty} \lec \norm{\de z}_{X_4}$.
\end{enumerate}
\end{remark}

\subsection{Small energy case}
\label{lwp:S4}

The proof of Theorem~\ref{lwp:T3}
is similar to that of Theorem \ref{lwp:T1}.

\begin{proof}
When $m \ge 1$, the limits $\lim_{ r \to 0} \bv_0(r)$ and $\lim_{ r \to
\infty} \bv_0(r)$ exist and it is necessary that $\bu_0(0) =
\bu_0(\infty)$.
We may assume $\bu_0(0) = \bu_0(\infty)=- \hk$.  In the
proof for Theorem \ref{lwp:T1}, we may redefine
\[
  \bh(r):= - \hk, \quad \bv(r) = (z_2, z_1, - 1 - \gamma)^T,
\]
and the parameters $s$ and $\al$ are no longer needed. The same proof,
in particular the difference estimate $\nrm{\de \tilde q^\sharp
}_{Str[I]} \lec (\si^{1/2}+ \de) \nrm{\de \tilde q}_{Str[I]}$, then
gives the local wellposedness.
\end{proof}

Note that this proof does not directly apply to the radial case, since
$\norm{\bu(r)}_{\dot H^1}$ no longer controls $\norm{z/r}_2$.


\section{Appendix: some lemmas}
\label{sec:lemmas}


\subsection{Computation of nonlinear terms}
\label{sec:nonlin}

To find the equations for $\dot s$ and $\dot \al$, we need to
compute $(g,({\bf V}^{\bh})^{-1}{P^{\bh} \bf F_1})_{L^2}$ for
$g=h_1$ or $g=N_0h_1$. Here is the result.

\begin{lemma} \label{Ap2:T2}
Recall ${\bf F_1} = - 2 \ga_r \frac mr h_1 \hj + \xi \times
(\Delta_r + \frac {m^2}{r^2} R^2 \xi)$ and $ ({\bf
V}^{\bh})^{-1}{P^{\bh} \bf F_1} = \hj \cdot {\bf F_1} + i ( h
\times \hj ) \cdot {\bf F_1}$. For any suitable function $g$,
\begin{align} \nonumber
(g,({\bf V}^{\bh})^{-1}{P^{\bh} \bf F_1})_{L^2}
 = \int_0^\infty \bigg( & i g_r (- \gamma z_r + z \gamma_r)
+ \frac mrh_1 g (-2 \gamma_r - i z_2 z_{1,r} + i z_1 z_{2,r})
\\
& + \frac mr(h_1g)_r (\gamma^2 - i z_2 z)
+ i \frac {m^2}{r^2}(2h_1^2-1) g \gamma z \bigg) rdr.
\end{align}
\end{lemma}

\begin{proof}
Decompose
\[
\int_0^\infty g ({\bf V}^{\bh})^{-1}{P^{\bh} \bf F_1} rdr = \int -
2 g \frac mr h_1 \ga_r + \int gP(\xi \times \Delta_r \xi) + \int
gP(\xi \times \frac {m^2}{r^2}R^2 \xi) = : I_1+I_2+I_3.
\]
Denote $[a,b,c]= a \hj + b h \times \hj + c h$.
For any vector $\eta$,
\[
P(\xi \times \eta) = [1,i,0] \cdot ([z_1,z_2,\ga] \times \eta)
=  ([1,i,0] \times  [z_1,z_2,\ga]) \cdot \eta
= [i \ga, - \ga, -iz] \cdot \eta.
\]
Since $h_r = \frac mr h_1 h \times \hj $,
\[
\pd_r [a,b,c] = [a_r, \ b_r+\frac mr h_1 c, \ c_r - \frac mr h_1 b].
\]
Thus
\begin{align*}
I_2 &= \int g[i\ga,-\ga,-iz] \cdot \Delta_r [z_1,z_2,\ga]
\\
&= \int \pd_r[-ig\ga,g\ga,igz] \cdot \pd_r [z_1,z_2,\ga]
\\
&= \int [-i (g\ga)_r, (g\ga)_r + \frac mr h_1 igz,\
 i(gz)_r - \frac mr h_1 g \ga]
\cdot [z_{1,r}, \ z_{2,r} + \frac mr h_1 \ga, \ \ga_r - \frac mr h_1 z_2]
\\
&= \int g  (-i \ga_r z_{1,r} + \ga_r z_{2,r} + iz_r \ga_r)
+ \int g_r (-i \ga   z_{1,r} + \ga   z_{2,r} + iz   \ga_r)
\\
& \quad + \int g_r(\frac mr h_1 \ga^2 - i \frac mr h_1 z_2 z)
+ \int g( - i \frac mr h_1 z_2 z_r + i \frac mr h_1 z z_{2,r})
+ \int g \frac {m^2}{r^2} h_1^2iz_1 \ga .
\end{align*}
Note that the first integral is zero, and we have canceled two
$\int g \frac mr h_1 \ga \ga_r$. Also,
\begin{align*}
I_3 &=  \int g[i \ga, - \ga, -iz] \cdot \frac {m^2}{r^2} R^2\xi =
\int g(\ga h_3-ih_1z,\ i\ga,\ *) \cdot  \frac {m^2}{r^2}
 (z_2h_3 - \ga h_1,\ z_1,\ 0)
\\
&=\int \frac  {m^2}{r^2} g( h_3^2 \ga z_2 - h_1h_3\ga^2
+ i h_1h_3 z_2 z + i h_1^2\ga z - i\ga z_1).
\end{align*}
Summing up $I_1+I_2+I_3$, we get the Lemma.
\end{proof}


\subsection{Linear weighted-$L^2$ estimate}
\label{sec:weighted}

\begin{lemma}
Let $H$ be a self-adjoint operator on $L^2(\R^n)$ satisfying
the weighted resolvent estimate
\[
  \sup_{\mu \not\in \R ; \;\; \phi \in L^2, \| \phi \|_{L^2} = 1}
  \norm{\frac{1}{|x|} (H - \mu)^{-1} \frac{1}{|x|} \phi}_{L^2} \lec 1.
\]
Then for $f(x,t) \in \frac{1}{|x|} L^2$,
\[
  \norm{ \frac{1}{|x|} \int_0^t e^{i(t-s)H} f(x,s) ds
  }_{L^2_{x,t}(\R^n \times \R)} \lec
  \norm{ |x| f }_{L^2_{x,t}(\R^n \times \R)}.
\]
\end{lemma}

\begin{proof}
First some simplifications.
It suffices to prove the estimate for 
$f(x,t)$ compactly supported, and 
$f(x,t) \in \frac{1}{|x|}L^2_{x,t} \cap L^\infty_t L^2_x$ 
(by density).
Also, it is enough to consider $t \geq 0$ (i.e. $f(x,t)$ 
supported in $\{t \geq 0\}$). 
Finally, we regularize the integral: set
\[
  F_\epsilon(x,t) := \frac{1}{|x|} \int_0^t
  e^{i(t-s)(H + i\epsilon)} f(x,s) ds.
\]
We will prove the estimate for $F_\epsilon$ with an
$\epsilon$-independent constant, and the lemma follows from this.
Under our assumptions, $F_\epsilon$ is well-defined
as a $\frac{1}{|x|} L^2_x$-valued function of $t$, and
\[
  \int_0^\infty \| |x| F_\epsilon(\cdot,t) \|_{L^2_x} dt < \infty.
\]
Hence the Fourier transform of $F_\epsilon$
in $t$ is well-defined (as a $\frac{1}{|x|} L^2_x$
-valued function of $\tau$):
\[
  \tilde{F}_\epsilon(x,\tau) := (2\pi)^{-1/2}
  \int_0^\infty e^{-i t \tau} F_\epsilon(x,t) dt.
\]
Changing the order of integration, we see
\[
\begin{split}
  \tilde{F}_\epsilon(x,\tau) &= \frac{1}{|x|} (2\pi)^{-1/2}
  \int_0^\infty dt e^{-i t \tau} \int_0^t ds
  e^{i(t-s)(H + i\epsilon)} f(x,s) \\ &=
  \frac{1}{|x|} (2 \pi)^{-1/2} \int_0^\infty ds  
  e^{-is(H + i \epsilon)}
  \int_s^\infty dt e^{it(H - \tau + i\epsilon)} f(x,s) \\
  &= \frac{1}{|x|} (2\pi)^{-1/2} (i) (H - \tau + i \epsilon)^{-1} 
  \int_0^\infty ds e^{-is\tau} f(x,s) ds \\
  &= \frac{1}{|x|}(i) (H - \tau + i \epsilon)^{-1} \tilde{f}(x,\tau)
\end{split}
\]
and so using the weighted resolvent estimate gives
\[
  \| \tilde{F}_\epsilon \|_{L^2_x} \lec
  \| |x| \tilde{f}(x,\tau) \|_{L^2_x},
\]
and squaring and integrating in $\tau$ yields
\[
  \| \tilde{F}_\epsilon \|_{L^2_{x,\tau}}^2 \lec
  \| |x| \tilde{f} \|_{L^2_{x,\tau}}^2 
  \lec \| |x| f \|_{L^2_{x,t}}^2.
\]
By a vector-valued version of the Plancherel theorem
(see eg. \cite{RS4}, Ch. XIII.7 ),
\[
  \| F_{\epsilon} \|_{L^2_{x,t}}^2 =
  \| \tilde{F}_\epsilon \|_{L^2_{x,t}}^2
  \lec \| |x| f \|_{L^2_{x,t}}^2,
\]
completing the proof.
\end{proof}


\subsection{Proof of the splitting lemma}
\label{sec:split}

Here we prove Lemma~\ref{lem:split}.

\begin{proof}
For $\bu = e^{m \theta R} \bv(r) \in \Sigma_m$, $s > 0$,
and $\alpha \in \R$, define
\[
  F(\bu;s,\alpha) := \int_0^\infty
  (\hj + i J^{\bh(\rho)}\hj) \cdot
  e^{-\alpha R} \bv(s\rho) h_1(\rho) \rho d \rho \;\; \in \C.
\]
Note that for $\bu$ of the form~\eqref{eq:uform},
\eqref{eq:orth2} is equivalent to $F(\bu;s,\alpha) = 0$.

Suppose $\sE(\bu) \le 4\pi m + \delta^2$.
It is shown in~\cite{GKT} that if $\delta$ is sufficiently small,
then there are $\hat{s}$, $\hat{\alpha}$, and $\hat{z}$ such that
$\bu(r,\theta) = e^{[m\theta + \hat{\alpha}]R}
[(1 + \hat{\gamma}(r/\hat{s}))\bh(r/\hat{s})
+ \bV^{r/\hat{s}}(\hat{z}(r/\hat{s}))]$, and with
$\| \hat{z} \|_X^2 \lec \delta_1^2 := \sE(\bu(0)) - 4\pi m
\le \delta^2$
(but $\hat{z}$ does not satisfy~\eqref{eq:orth2}).

It follows from this, and the fact that $\rho h_1(\rho) \in L^2(\rho d\rho)$
for $m \geq 3$, that for some $\delta_0 > 0$, $F$ is a $C^1$ map from
\[
  \{ \bu \in \Sigma_m \; | \; \sE(\bu) \leq 4\pi m + \delta_0^2 \}
  \times (\R^+ \times \R)
\]
into $\C$.  Furthermore, straightforward computations show that
\[
  F(e^{m\theta R} \bh(r); 1, 0) = 0,
\]
and
\[
  \left( \begin{array}{c}
  \partial_s F(e^{m\theta R} \bh(r); 1, 0) \\
  \partial_{\alpha} F(e^{m\theta R} \bh(r); 1, 0)
  \end{array} \right) =
  \| h_1 \|_{L^2}^2 \left( \begin{array}{c} i \\ -1
  \end{array} \right).
\]
By the implicit function theorem, we can solve
$F = 0$ to get $s = s(\bu)$ and $\alpha = \alpha(\bu)$ for $\bu$ in a
$\dot H^1$-neighbourhood of the harmonic map $e^{m \theta R}\bh(r)$.

Since $\| \hat{z} \|_X \lec \delta$,
provided $\delta$ is chosen small enough (depending on the
size of this neighbourhood),
\[
  \hat{\bu}(x) := e^{-\hat{\alpha}R} \bu(\hat{s} x)
  = e^{m\theta R}[(1 + \hat{\gamma}(r)) \bh(r) + \bV^{r} \hat{z}(r)]
\]
lies in this neighbourhood, yielding $s(\hat{\bu})$
and $\alpha(\hat{\bu})$ with
$F(\hat{\bu}; s(\hat{\bu}), \alpha(\hat{\bu})) = 0$.
Furthermore,
\[
  |s(\hat{\bu}) - 1| + |\alpha(\hat{\bu})|
  \lec \| \hat{z} \|_X,
\]
and so
\[
  \| z(\rho) \|_X =
  \| (\hj + i J^{\bh(\rho)} \hj) \cdot e^{-\alpha(\hat{\bu}) R}
  \hat{\bv}(s(\hat{\bu}) \rho) \|_X
  \lec \| \hat{z} \|_X \lec \sE(\bu) - 4\pi m.
\]

To complete the proof of the lemma, undo the rescaling: set
$s(\bu) := s(\hat{\bu})/\hat{s}$ and
$\alpha(\bu) := \alpha(\hat{\bu}) + \hat{\alpha}$.
\end{proof}


\donothing{
Let $\frac 1r L^2_{rad}$ denote the space of all functions $ f(r)$
with $rf(r)\in L^2_{rad}$.

\begin{lemma}[The range of $N_0$] Let $m \ge 1$ be any positive integer and
$N_0 = - \Delta _r + \frac {m^2}{r^2}$.

(i) For all $ f \in \frac 1r L^2_{rad}$, there exists $g \in X$ so
that $N_0 g = f$.

(ii) Consider $N_0$ as an operator on $X$. Then $X \cap \frac 1r L^2_{rad}
\subset Ran(N_0)$.

\end{lemma}

\begin{proof}
It suffices to show (i) since it implies (ii).
Consider real-valued functions on $(0,\infty)$.  The general solution
of the ODE $N_0g = 0$ is $g(r)=a r^m + b r^{-m}$.  The ODE $N_0g =f $
for any function $f(r)$ can be solved by the method of variations of
parameters:
\[
g(r) = r^m a(r) + r^{-m} b(r), \quad
g'(r) = \frac mr [r^m a(r) - r^{-m} b(r)],
\]
where
\[
a(r)=- \frac 1{2m}\int_{r_1}^r s^{1-m} f(s) ds , \quad
b(r)=\frac 1{2m}\int_{r_2}^r s^{1+m} f(s) ds ,
\]
for some constants $r_1,r_2$. We will choose $r_1 = \infty$ and $r_2 = 0$.
The condition $g \in X$, i.e., $g',g/r \in L^2_{rad}$, is equivalent to
\begin{equation}\label{eq1}
r^{m-1}\int_{\infty}^r s^{1-m} f(s) ds\in L^2_{rad} \quad \text{and} \quad
r^{-m-1}\int_{0}^r s^{1+m} f(s) ds  \in L^2_{rad}.
\end{equation}
Let $F(r) =\int_{\infty}^r s^{1-m} f(s) ds$, which vanishes as $ r \to
\infty$. (We may assume $f \in C^1_c(0,\infty)$ by density.) The first
statement of \eqref{eq1} is true if
\[
\|r^{m-1}F\|_{L^2_{rad}} \lec  \|r^{m}F'\|_{L^2_{rad}} =\| rf\|_{L^2_{rad}} ,
\]
but this follows from Hardy's inequality in $\R^{2m+2}$.  The second
statement of \eqref{eq1} follows from the first by duality: For any
$k(r) \in L^2_{rad}$,
\[
\int_0^\infty k(r) r^{-m-1}\int_{0}^r s^{1+m} f(s) ds \,r dr
=\int_0^\infty  sf(s)
\bke{s^{m-1} \int_{s}^\infty r^{1-m} \frac{k(r)}r  dr} s ds ,
\]
which is bounded by $C \|sf(s)\|_{L^2_{rad}} \cdot  \|k(r)\|_{L^2_{rad}}$.
\end{proof}
}

\section*{Acknowledgments}

The first and third authors are grateful for support from NSERC
grants. The second author was supported partly by the 
Korea Research Foundation Grant funded by the Korean 
Government (MOEHRD, Basic Research Promotion Fund)
(KRF-2006-003-C00020). 
Part of this work was completed while the second author
was supported by a PIMS postdoctoral fellowship at UBC.


Stephen Gustafson

Mathematics Department

University of British Columbia

Vancouver, BC, Canada, V6T 1Z2

gustaf@math.ubc.ca

\vskip.5cm

Kyungkeun Kang

Department of Mathematics

Sungkyunkwan University and Institute of Basic Science

Suwon 440-746, Republic of Korea

kkang@skku.edu

\vskip.5cm

Tai-Peng Tsai

Mathematics Department

University of British Columbia

Vancouver, BC, Canada, V6T 1Z2

ttsai@math.ubc.ca

\end{document}